\documentclass[11pt]{article}
\usepackage{amsfonts,color}
\usepackage{amssymb}
\usepackage{amsmath}
\usepackage{amsthm}
\usepackage{natbib}
\usepackage{enumerate}
\usepackage{graphicx}
\usepackage{wrapfig}
\usepackage{rotating}
\usepackage{epsfig}
\usepackage[left=1in, right=1in, top=0.8in, bottom=1in]{geometry}
\usepackage{mdwlist}
\usepackage{secdot}
\usepackage{comment}
\usepackage{mathrsfs}
\usepackage{mdwlist}
\usepackage{multirow}
\usepackage{lscape}
\usepackage{pdflscape}
\usepackage{array}
\usepackage[compact]{titlesec}
\usepackage{times}
\usepackage{tikz}
\usepackage{subfigure}
\usepackage[justification=centering]{caption}
\usepackage{setspace}
\usepackage{color}
\usepackage{authblk}
\usepackage{wrapfig}
\usepackage{hyperref}
\usepackage{makeidx}

\bibpunct{[}{]}{,}{n}{}{;}
\titlespacing*{\section}{0pt}{*2}{*0}
\titlespacing*{\subsection}{0pt}{*2}{*0}
\titlespacing*{\subsubsection}{0pt}{*1}{*0}

\newtheorem{theorem}{Theorem}
\newtheorem{lemma}{Lemma}

\newtheorem{proposition}{Proposition}
\newtheorem{corollary}{Corollary}

\newtheorem{definition}{Definition}

\newtheoremstyle{compacttheorem}
{3pt}
{3pt}
{}
{}
{}
{}
{}
{}

\title{The Covering Path Problem on a Grid}
\author[1]{\rm Liwei Zeng}
\author[2]{\rm Sunil Chopra}
\author[1]{\rm Karen Smilowitz}
\affil[1]{Department of Industrial Engineering and Management Sciences}
\affil[2]{Kellogg School of Management}
\affil[ ]{Northwestern University}

\begin{document}
\maketitle
\begin{abstract}
This paper introduces the covering path problem on a grid (CPPG) which finds the cost-minimizing path connecting a subset of points in a grid such that each point that needs to be covered is within a predetermined distance of a point from the chosen subset. We leverage the geometric properties of the grid graph which captures the road network structure in many transportation problems, including our motivating setting of school bus routing. As defined in this paper, the CPPG is a bi-objective optimization problem comprised of one cost term related to path length and one cost term related to stop count. We develop a trade-off constraint which quantifies the trade-off between path length and stop count and provides a lower bound for the bi-objective optimization problem. We introduce simple construction techniques to provide feasible paths that match the lower bound within a constant factor. Importantly, this solution approach uses transformations of the general CPPG to either a discrete CPPG or continuous CPPG based on the value of the coverage radius. For both the discrete and continuous versions, we provide fast constant-factor approximations, thus solving the general CPPG.
\end{abstract}

\textbf{Keywords:} covering path problem; grid optimization; school bus routing; location routing problem

\section{Introduction} \label{introduction}
School bus routing is an activity performed by school districts across the United States, with annual costs over \$20 billion (\cite{national_costs}). The core subproblem in school bus routing considers an area where children are to be picked up by a bus. The goal is to identify \textit{bus stops} such that each child is sufficiently close to a stop and a \textit{bus route} such that the total cost/time of the bus is minimized. Both travel and stops incur cost and time for the bus. Our goal in this paper is to find high quality solutions to the sub-problem quickly. This will allow decision makers to interactively change parameters, such as the area assigned to a bus or the definition of ``sufficiently close'',  to evaluate the extent to which such changes impact total cost.

The covering path problem (CPP) that has been studied in the literature captures many elements of the core sub-problem in school bus routing.  The CPP is a variant of the traveling salesman problem (TSP), in which the vehicle is not required to visit every point and the path does not end at the starting point. Like the TSP, the CPP is NP-hard on general graphs (\cite{current1981multiobjective}). We develop efficient solution methods for the CPP when the problem is restricted to a grid graph. This is an important restriction that naturally arises when routing school buses in many urban and suburban settings where the underlying road network resembles a grid.

\begin{definition}[\textbf{CPP}]\label{def:CPP}
	Consider a graph $G=(V, E)$ with edge weights $l_{e}$ for $e\in E$ and node weights $t_{v}$ for $v\in V$, a coverage region $\mathcal{R}$, and a coverage radius $k>0$. The CPP finds a set of stops $V_{1}\subseteq V$ such that for every point $x\in \mathcal{R}$, there exists $v_{1}\in V_{1}$ at a distance at most $k$, and the minimum cost path $P$ connecting the nodes in $V_{1}$. Given scalars $L\in \mathbb{R}$ and $T\in \mathbb{R}$ and any function $C(L, T)$, the cost of path $P$ is given by
	\begin{equation}
		Cost(P)=C(L, T)=C(\sum_{e\in P}l_{e}, \sum_{v\in V_{1}}t_{v}).
	\end{equation}
	
\end{definition}

A point $v_{1}\in V_{1}$ is referred to as a \emph{stop} and a path connecting all nodes in $V_{1}$ is referred to as a \emph{covering path}. Point $A$ is said to \emph{cover} point $B$ if and only if the distance between $A$ and $B$ is no more than $k$. The two cost terms $L$ and $T$ are referred to as \emph{path length} and \emph{stop count}, respectively.  In the CPP literature, the coverage region $\mathcal{R}$ is typically a set of nodes to be covered (which may be the set $V$) and the stops $V_{1}$ are chosen from that node set.  In our work for school bus routing, we interpret the coverage region as the area in which students live and choose bus stops from a discrete subset of nodes in the area.

This paper is motivated by a collaboration with a public school district focused on improving service and lowering cost for bus transportation. The underlying road network for the district resembles a grid and our goal is to leverage this structure to obtain robust transportation solutions, thus allowing the school district to (1) employ simple strategies to identify bus stop locations and plan bus routes and (2) embed these strategies and associated cost approximations in broader decision frameworks, covering decisions such as student assignments to schools.

The School Bus Routing Problem (SBRP) has been extensively studied in the literature (e.g., \cite{desrosiers1980overview}, \cite{newton1969design}, \cite{park2010school}). The SBRP itself is a composite of five subproblems: data preparation, bus stop selection, bus route generation, school bell time adjustment, and route scheduling. As noted earlier, students are often not picked up at their homes, but rather are assigned to bus stops within a set walking distance, thus the SBRP is studied as a combination of routing and covering. Existing SBRP literature typically does not specify the underlying graph structure of the road network. The joint subproblem of bus stop selection and bus route generation is modeled with integer programming models (\cite{bowerman1995multi}, \cite{desrosiers1980overview}, \cite{gavish1979approach}, \cite{park2010school}). As the joint subproblem is NP-Hard, finding solutions can be challenging for implementation in practice. Later work designs heuristics to tackle the computational complexity, such as genetic algorithm (\cite{diaz2012vertical}), tabu search (\cite{pacheco2013bi}) and randomized adaptive search procedure (\cite{schittekat2013}). Heuristics for the combined problem of bus stop selection and route generation mainly follow two strategies: the location-allocation-routing (LAR) strategy (\cite{bodin1979routing}) and the allocation-routing-location (ARL) strategy (\cite{bowerman1995multi}). The LAR strategy sequentially selects bus stops, assigns students to bus stops, and designs bus routes. The ARL strategy first groups students into clusters, selects stops and generates routes for each cluster, and then assigns students to bus stops. These strategies solve routing and location problems sequentially. In our work, we aim to solve these two subproblems simultaneously by leveraging the grid structure of the underlying graph. This study of the CPP in a stylized grid setting is a first step in our analysis of the joint problem of route design and bus stop selection. Our solution approach is particularly useful in a setting where the decision maker cares about both the number of stops and route length and wants to interactively adjust the weight assigned to either when designing routes. Our approach quickly provides a high quality solution to the decision maker. In the conclusion, we discuss next steps to use these results in a stylized setting to address more complex settings with multiple vehicles and other generalizations.

Motivated by the SBRP in which students are typically located along streets and bus stop locations are selected from intersections, we define the following notation. A unit grid graph, also known as a square grid graph (\cite{weisstein2001grid}), is a graph whose nodes correspond to integer points in the plane with the $x$-coordinates from 0 to $m$ and $y$-coordinates from 0 to $n$. Two nodes in the grid graph are connected if and only if they are end points of an edge of distance 1. Given that students live on streets that correspond to edges of the grid graph and walk primarily along the streets, we use the $l_1$ norm to measure the distance between any two points. The $l_1$ norm is the shortest path length between two points traveling only along edges in the grid graph. We show in our analysis that the grid structure leads to strong approximation results (in some cases near-optimal) for the optimization problem with a linear objective function of the number of stops and the route length. We define the unit grid formally below.

\begin{definition}[\textbf{Unit Grid}]
	Given $m,n\in \mathbb{N}^{+}$, let $D_{int}$ be the set of integer points $(x,y)$ with $0\leq x\leq m$ and $0\leq y\leq n$. Let $E$ be the set of grid edges connecting nodes $(x,y)$ in $D_{int}$ with adjacent grid nodes. $G=(D_{int}, E)$ defines a $m\times n$ unit grid.
\end{definition}

The vertices of the unit grid, $V$, are equivalent to the set $D_{int}$, the set of integer points in the grid. In a unit grid, $l_{e}=1$ for $e\in E$ and $t_{v}=1$ for $v\in V$. We further define $D_{edge}$ to be all points on the edges $E$ of the grid. Observe that $D_{edge}$ includes not only the vertices $V$ but also all points on the edges $E$. CPPG can now be specialized from the definition of CPP as follows:

\textbf{CPPG (Covering Path Problem on a Grid)}.  Solve CPP given  $m\times n$ unit grid $G=(V,E)$ with $\mathcal{R}=D_{edge}$ and $V_1 \subseteq D_{int}$.

In practice, it may be hard for decision makers to provide an exact cost function $C(L,T)$. Thus, we look to solve the optimization problem of minimizing $C(L,T)$ without specifying the form of the function $C(\cdot,\cdot)$. In order to do so, we consider a related decision version of the CPPG. By fully characterizing the decision problem, we can quantify the trade-off between the path length $L$ and the number of stops $T$, thereby efficiently solving the bi-objective optimization problem.

\textbf{Decision version of CPPG}.  Given $L, T\in \mathbb{R}$, determine if there exists a set of stops $V_{1}$ and a covering path $P$ such that $\sum_{e\in P}l_{e}\leq L$ and $\sum_{v\in V_{1}}t_{v}\leq T$.

\begin{wrapfigure}{r}{0.5\textwidth}
	\begin{center}
		\includegraphics[height=2in]{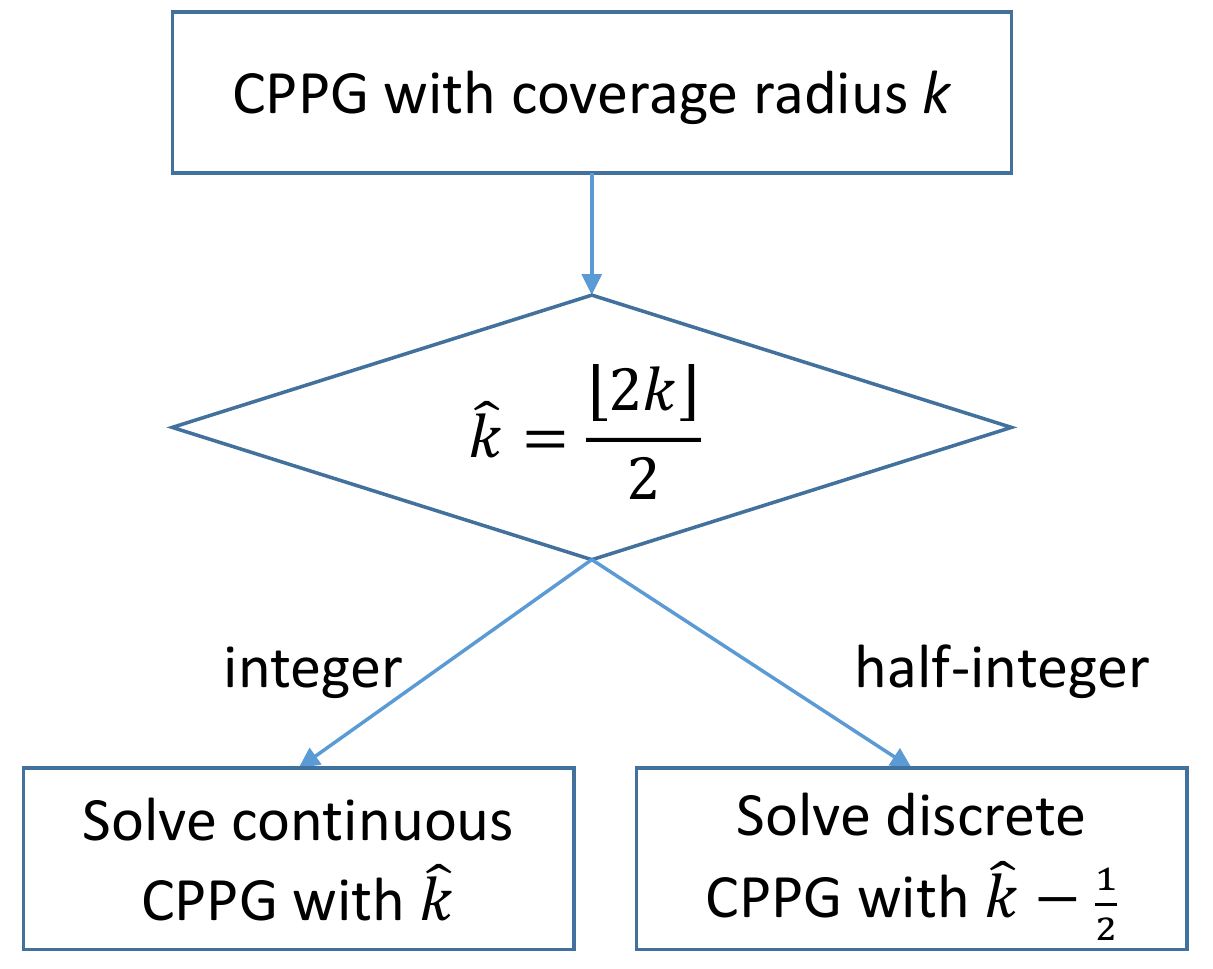}
		\caption{Solution approach for CPPG}
		\label{fig:solapp}
	\end{center}
\end{wrapfigure}

To solve the CPPG, we first show that one can reduce the CPPG with any coverage radius $k > 0$ to one of two cases: $k$ as integer or half-integer. Recall that the coverage region in the CPPG is the set of points on the edges $D_{edge}$. We show that: (1) when $k$ is an integer, covering $D_{edge}$ is equivalent to covering the rectangle $D=\{(x,y), 0\leq x\leq m, 0\leq y\leq n\}$ with coverage radius $k$; (2) when $k$ is a half-integer, covering $D_{edge}$ is equivalent to covering integer points $D_{int}$ with coverage radius $k-\frac{1}{2}$. The transformation then leads to two variations of the CPPG: the continuous CPPG in which we cover all points in the rectangular grid $D$ and the discrete CPPG in which we cover all integer points in $D_{int}$ (see Figure \ref{fig:solapp} for illustration). The continuous CPPG falls into a stream of continuous facility location and routing models which have been shown to offer computational simplicity compared to their discrete counterparts. We show that insights from the continuous CPPG can be used for the discrete CPPG which can be viewed as a CPP on a grid graph in which only a finite number of points must be covered. We develop efficient methods to find feasible, high quality solutions for both variations, thus solving the original problem.

We solve both covering path problem variants in a bi-objective setting: i.e., we minimize a function of the path length $L$ and stop count $T$. Similar covering tour problems have been studied in the literature: \cite{jozefowiez2007bi} seek to minimize two objectives--tour length and coverage radius and \cite{tricoire2012bi} study the trade-off between fixed cost and uncovered demand in a stochastic setting by characterizing the Pareto frontier. We solve the bi-objective problems by identifying the set of all non-dominated solutions (or Pareto frontier) with an inequality that quantifies the trade-off between path length and stop count. Our approach allows us to find high quality solutions quickly which is necessary in an interactive setting where the decision maker evaluates different weights on length and number of stops. By high quality, we mean solutions that are within a fixed ratio of the optimal solution. By quickly, we mean polynomial time.

The remainder of this paper is organized as follows.  In Section \ref{Literature review} we review related work on the CPP and optimization problems on grid graphs. In Section \ref{initial} we establish the transformations needed for the approach in Figure \ref{fig:solapp} and formally introduce the continuous and discrete CPPG. In Section \ref{relax} we present a relaxation of the continuous CPPG that leads to foundational results which are used in Sections \ref{continuous} and \ref{discrete} for the continuous and discrete CPPG, respectively. Finally, we conclude in Section \ref{conclusion} with next steps to apply these results to the SBRP.

\section{Literature Review} \label{Literature review}
We review two relevant streams of related research: covering tour and path problems and optimization problems on grid graphs.

\subsection{Covering tour and path problems}
\cite{current1981multiobjective} introduces the CPP and shows its NP-hardness from a reduction of the TSP when the coverage radius equals to zero. The covering tour problem (CTP) is similar to the CPP, requiring the path to start and end at the same point. The CPP can be reduced to the CTP by adding a dummy node that is connected to all other nodes with zero cost but not covered by any other nodes. Existing work formulates the CTP as an integer linear program (ILP), beginning with \cite{current1989covering} and builds corresponding solution approaches. \cite{gendreau1997covering} study the polyhedron of the ILP and provide a branch-and-cut algorithm. \cite{hachicha2000heuristics} present an ILP formulation and heuristics for the multi-vehicle CTP. The CTP can also be treated as a generalized traveling salesman problem (GTSP) (\cite{fischetti1997branch}): given several sets of nodes, the GTSP seeks to determine a shortest tour passing at least once through each set. Recent work continues on designing solution approaches for the CTP and multi-vehicle CTP; such as branch and price (\cite{jozefowiez2014branch}), column generation (\cite{murakami2014column}) and adaptive large neighborhood search (\cite{leticia2015selector}). The CTP has also been studied in the bi-objective setting. \cite{jozefowiez2007bi} introduce the bi-objective CTP which aims to minimize both the tour length and the coverage radius. \cite{tricoire2012bi} study the stochastic bi-objective CTP and discuss the fundamental trade-off between fixed cost and uncovered demand. Different from previous work, we develop polynomial solution methods that exploit the underlying grid structure to obtain provable bounds.

Combining facility location and route design has also been broadly studied in other related problems. The location-routing problem (LRP) pays special attention to the underlying issue of vehicle routing (see \cite{albareda2015location}, \cite{drexl2015survey}, \cite{prodhon2014survey} for reviews). There have been several different formulations of the capacitated LRP in recent work (\cite{cherkesly2017set}, \cite{contardo2013computational}, \cite{contardo2013exact}) where the coverage distance is relaxed and the distance to nodes not on a path becomes a cost to minimize. In the ringstar problem (\cite{labbe2004ring}), the objective function combines location and routing costs with  a linear combination of path cost and access cost from assigning nodes to facilities.

Recently, attention has been given to continuous facility location problems with access costs. \cite{carlsson2014continuous} introduce a problem that is related to our work. They consider facility location with backbone network costs, where the objective function is a linear combination of fixed costs from installing facilities, backbone network costs from connecting facilities and access costs from connecting customers to facilities.  The fixed costs in \cite{carlsson2014continuous} are equivalent to our fixed costs of stops and the backbone network is equivalent to our covering path. In our problem, the access cost is modeled as a coverage constraint where each stop covers points within a given distance, consistent with how the problem is viewed by the school district. Our use of the $l_1$ norm to calculate distance is also consistent with the practical problem. In this setting, we develop a solution approach that provides high quality approximation solutions when the objective function is increasing and convex. Our solution approach characterizes the boundary of all feasible solutions and uses this characterization to solve the optimization problem.

\subsection{Optimization on grid graphs}
Some of the most fundamental combinatorial optimization problems have been well studied on grid graphs and other graphs with special metrics and topological structures. For the Hamiltonian Cycle Problem, \cite{itai1982hamilton} prove its NP-hardness on general grid graphs and \cite{umans1997hamiltonian} show the problem can be solved in polynomial time on a simple grid graph (e.g., a grid without holes). For the TSP, \cite{arora1998polynomial} provides a polynomial-time approximation scheme for problems on grid graphs but the algorithm is computationally challenging for large instances. Recent advances in the TSP also indicate potential benefits of working on structured graphs. \cite{arkin2000approximation} give a $\frac{6}{5}$-approximation polynomial-time algorithm for the TSP on a simple grid graph. \cite{gharan2011randomized} provide a $(\frac{3}{2}-\varepsilon)$-approximation polynomial-time algorithm for the graph TSP where edge cost is measured by the shortest path length on a unit graph. Their algorithm follows the structure of Christofides heuristic by cleverly choosing a random spanning tree (not always the minimum spanning tree) in the first step. The results improve the $\frac{3}{2}$-approximation due to \cite{christofides1976worst} for this TSP variant. \cite{sebo2014shorter} later improve the approximation ratio to $\frac{7}{5}$ together with a derandomized algorithm using forest representations of hypergraphs. \cite{savacs2002finite} study the facility location problem with barriers using the $l_1$ norm to measure distance. In general, the grid assumption provides a unified geometric structure with fewer degrees of freedom. Such structures can be easier to analyze with the help of geometric and combinatorial techniques. Motivated by these results, our paper looks to leverage the grid structure to solve the CPP.

Our contribution in this paper is to provide a polynomial algorithm for the CPPG that exploits the underlying grid structure to obtain high quality solutions for an objective function that accounts for both the number of stops and the path length. This bi-objective scenario naturally arises as the core problem in school bus routing.

\section{Characterizing CPPG Problem Settings} \label{initial}
In this section, we present preliminaries for our CPPG solution approach in Figure \ref{fig:solapp}. On a unit grid graph, when $k<1$, CPPG has a trivial solution where all nodes in V are stops, i.e., $V_1=V$. Thus, for the rest of the paper we assume that $k\geq 1$. We show that one can round the coverage radius $k$ down to the nearest integer or half-integer and maintain the coverage properties. When this rounding results in an integer value, the CPPG can be solved with the continuous CPPG and when the rounding results in a half-integer value, the CPPG can be solved with the discrete CPPG.

\subsection{Rounding the coverage radius}
Recall the rectangular region $D=\{(x,y)~|~0\leq x\leq n, 0\leq y\leq m \}$ in $\mathbb{R}^{2}$ with a grid graph with the set of points on the edges $D_{edge}=\{(x,y)~|~(x,y)\in D, x\in \mathbb{Z}~\textrm{or}~y\in \mathbb{Z}\}$ and integer point set $D_{int}=\{(x,y)~|~0\leq x\leq n, 0\leq y\leq m, x,y\in \mathbb{Z}\}$. Without loss of generality, we assume  $m\geq n>0$. The following proposition shows that to solve the CPPG, it is sufficient to model $k$ as either integer or half-integer.

\begin{proposition}\label{prop:1}
	In the CPPG, a covering path with coverage radius $k$ is also a covering path with coverage radius $\frac{\left\lfloor 2k \right\rfloor}{2}$.
\end{proposition}

\emph{Proof of Proposition \ref{prop:1}}. Let $V_{1}\subseteq D_{int}$ be the set of stops in a covering path with coverage radius $k$. Given two points $x$ and $F$, let $||x-F||_1$ represent the $l_1$ distance between $x$ and $F$.

For each point $x\in D_{edge}$, let $$dist(x)=\textrm{min}_{F\in V_{1}}||x-F||_{1}$$ be the distance from $x$ to its nearest stop. All points in $D_{edge}$ are covered with radius $k$ if and only if  $\textrm{max}_{x\in D_{edge}}\{dist(x)\}\leq k$. We prove that the value of $\textrm{max}_{x\in D_{edge}}\{dist(x)\}$ is an integer or a half-integer. Therefore, $\textrm{max}_{x\in D_{edge}}\{dist(x)\}\leq \frac{\left\lfloor 2k \right\rfloor}{2}$.

Let $a,b\in D_{int}$ be two connected integer points in the grid graph. Given that all stops are located at integer points, $dist(a)$ and $dist(b)$ must be integers. Because $a$ and $b$ are connected, $|dist(a)-dist(b)|\leq ||a-b||_{1}=1$, which leads to three cases when we consider $\overrightarrow{ab}$, the edge connecting $a$ and $b$.

\emph{Case 1:} $dist(a)=dist(b)=c\in \mathbb{Z}$. For $x\in \overrightarrow{ab}$, the function $dist(x)$ takes maximum value of $c+\frac{1}{2}$ (which is a half-integer) at the midpoint of edge $\overrightarrow{ab}$.

\emph{Case 2:} $dist(a)-dist(b)=1$. For $x\in \overrightarrow{ab}$, the function $dist(x)$ is linear on edge $\overrightarrow{ab}$ and takes maximum value at point $a$.

\emph{Case 3:} $dist(a)-dist(b)=-1$. This case is symmetric to \emph{Case 2} and $dist(\cdot)$ takes maximum value at point $b$.

In summary, the maximum function value of $dist(\cdot)$ is either an integer or a half-integer on each edge. Therefore, the coverage property remains unchanged after rounding $k$ down to $\frac{\left\lfloor 2k \right\rfloor}{2}$.
\qed

Proposition \ref{prop:1} states that the coverage radius $k$ can be reduced to $\frac{\left\lfloor 2k \right\rfloor}{2}$, which is the largest integer or half-integer less than or equal to $k$. When $\frac{\left\lfloor 2k \right\rfloor}{2}$ is an integer, we expand the coverage region from $D_{edge}$ to the full rectangle $D$ so that the problem falls into the stream of continuous facility location; when $\frac{\left\lfloor 2k \right\rfloor}{2}$ is a half-integer, we restrict the coverage region to integer points $D_{int}$ and the problem can be viewed as a CPP on a grid graph in which we only cover a finite number of points. By replacing $D_{edge}$ with $D$ and $D_{int}$, we are able to obtain tighter approximation results for both the optimization and decision versions of the CPPG.

\subsection{Solving the CPPG with integer coverage radius}
As a first step, we show that when the coverage radius $k$ is an integer, solving the CPPG with coverage region $D_{edge}$ is equivalent to solving the CPPG with coverage region $D$.

\begin{proposition}\label{prop:2}
	Given a positive integer $k$, covering all points on the edges, $D_{edge}$, with radius $k$ is equivalent to covering all points in the rectangle, $D$, with the same radius.
\end{proposition}
\emph{Proof of Proposition \ref{prop:2}.}
Given that $D_{edge} \subset D$, covering all points in $D$ with radius $k$ naturally covers all points in $D_{edge}$ with radius $k$. To prove the other direction, it suffices to show that if $dist(x)\leq k~\forall x\in D_{edge}$, then $dist(x)\leq k~\forall x\in D$. We call $y\in D$ a mid-integer point if one of its coordinates is an integer and the other a half-integer; i.e., $y$ is the midpoint of two integer points with distance 1. Given $x\in D$, let $y_{x}\in D$ be the closest mid-integer point to $x$. We have $||x-y_{x}||_{1}\leq \frac{1}{2}$. Since $y_{x}$ is a mid-integer point, $dist(y_{x})$ must be a half-integer. Together with the fact that $y_{x}\in D_{edge}$, we have $dist(y_{x})\leq k-\frac{1}{2}$. From the triangle inequality, $dist(x)\leq dist(y_{x})+||x-y_{x}||_{1}\leq (k-\frac{1}{2})+\frac{1}{2}\leq k$. Thus, all points in the rectangle $D$ are covered with radius $k$.
\qed

With Proposition \ref{prop:2}, when $k$ is an integer, we solve the CPPG by defining an equivalent problem called the continuous CPPG, where the coverage region is expanded to the rectangle $D$.

\vspace{0.2in}


\textbf{C-CPPG (Continuous CPPG)}.  Solve CPP given an $m\times n$ unit grid $G=(V,E)$ with $\mathcal{R}=D$ and $V_1 \subseteq D_{int}$.

\vspace{0.2in}

\subsection{Solving the CPPG with half-integer coverage radius}
When the coverage radius $k$ is a half-integer, we show that solving the CPPG with coverage region $D_{edge}$ is equivalent to solving the CPPG with coverage region $D_{int}$, the set of integer points.

\begin{proposition}\label{prop:3}
	Given a half-integer $k$, covering all points in the edges, $D_{edge}$, with radius $k$ is equivalent to covering the integer points, $D_{int}$, with radius $k-\frac{1}{2}$.
\end{proposition}
\emph{Proof of Proposition \ref{prop:3}.}
For any integer point $z\in D_{int}$, $dist(z)$ is an integer less than or equal to $k$. Since $k$ is a half-integer, $dist(z)\leq k-\frac{1}{2}$. For any $x\in D_{edge}$, let $z_{x}$ be the nearest integer point to $x$. From the triangle inequality, $dist(x)\leq dist(z_{x})+||x-z_{x}||_{1}\leq (k-\frac{1}{2})+\frac{1}{2}=k$. Given that all $x\in D_{int}$ are within distance $k-\frac{1}{2}$ of stops, all $x\in D_{edge}$ must be within distance $k$ of stops.
\qed

With Proposition \ref{prop:3}, when $k$ is a half-integer, we solve the CPPG by defining an equivalent problem called the discrete CPPG, where the coverage region is restricted to the integer points $D_{int}$,  equivalent to the node set $V$.

\vspace{0.2in}

\textbf{D-CPPG (Discrete CPPG)}.  Solve CPP given an  $m\times n$ unit grid $G=(V,E)$ with $\mathcal{R}=D_{int}$ and $V_1 \subseteq D_{int}$.

\vspace{0.2in}

\subsection{A relaxation of the C-CPPG}
To develop solution approaches for the C-CPPG and D-CPPG, we define a relaxation of the C-CPPG, RC-CPPG,  where the set of stops $V_{1}$ can be selected from all points in $D$ (rather than $D_{int}$). The three key elements in the definition of each problem we study in the paper (potential stop locations, coverage region and coverage radius) are presented in Table \ref{tab:psettings}. Note in Table \ref{tab:psettings}, that choosing stop locations from the set of integers $D_{int}$ is equivalent to choosing from the node set $V$ in the unit grid.

\vspace{0.2in}


\textbf{RC-CPPG (A Relaxation of the C-CPPG)}.  Solve CPP given an $m\times n$ unit grid $G=(V,E)$ with $\mathcal{R}=D$ and $V_1 \subseteq D$.

\vspace{0.2in}

In Section \ref{relax}, we show that the optimization problem defined for the RC-CPPG can be solved to near-optimality. The continuity in the relaxed problem (both in potential stop locations and coverage region) allows for more accurate analysis and tighter bounds for the optimization and decision problems. The analysis approach developed for the RC-CPPG is applied to the C-CPPG in Section \ref{continuous} and D-CPPG in Section \ref{discrete} to obtain tight bounds.

\begin{table}[h]
	\begin{center}
		\begin{tabular}{ |c|c|c|c|c| }
			\hline
			Setting/Problem & CPPG  & C-CPPG & D-CPPG & RC-CPPG \\ \hline
			Potential stop locations & $D_{int}$  & $D_{int}$ & $D_{int}$ & $D$ \\ \hline
			Coverage region & $D_{edge}$ & $D$ & $D_{int}$ & $D$ \\ \hline
			Coverage radius & $k>0$ & $k\in \mathbb{N}^{+}$ & $k\in \mathbb{N}^{+}$ & $k>0$ \\ \hline
		\end{tabular}
		\caption{Differences in CPPG problem settings}
		\label{tab:psettings}
	\end{center}
\end{table}

\section{Analysis of the RC-CPPG \label{relax}}
The RC-CPPG falls into the regime of continuous facility location (\cite{daganzo1986configuration}). \cite{carlsson2014continuous} study a similar problem of finding a minimum cost path that covers a convex polygon. The objective function in \cite{carlsson2014continuous} is a linear combination of fixed cost, path cost and access cost. Both \cite{carlsson2014continuous} and our work handle the continuity of the coverage region with combinatorial and geometric approaches. While our paper only considers covering a rectangle with minimum fixed cost and path cost, we obtain stronger approximation results (in some settings near-optimal) by leveraging the structure of the $l_{1}$ norm and the grid graph.

The remainder of this section is organized as follows. We introduce the trade-off constraint which quantifies the trade-off between the path length $L$ and number of stops $T$ and provides lower bounds for the set of feasible solutions to the optimization problem. We then construct a family of feasible paths called ``up-and-down paths" that provide an upper bound that matches the lower bound obtained from the trade-off constraint. Finally, we present a polynomial algorithm that solves the optimization problem with a linear objective function in a near-optimal fashion. The proof techniques and results in this section are used in Sections \ref{continuous} and \ref{discrete} with slight changes.

\subsection{Trade-off constraint for the RC-CPPG}\label{sec:4.1}
\begin{wrapfigure}{r}{0.5\textwidth}
	\begin{center}
		\includegraphics[width=0.48\textwidth]{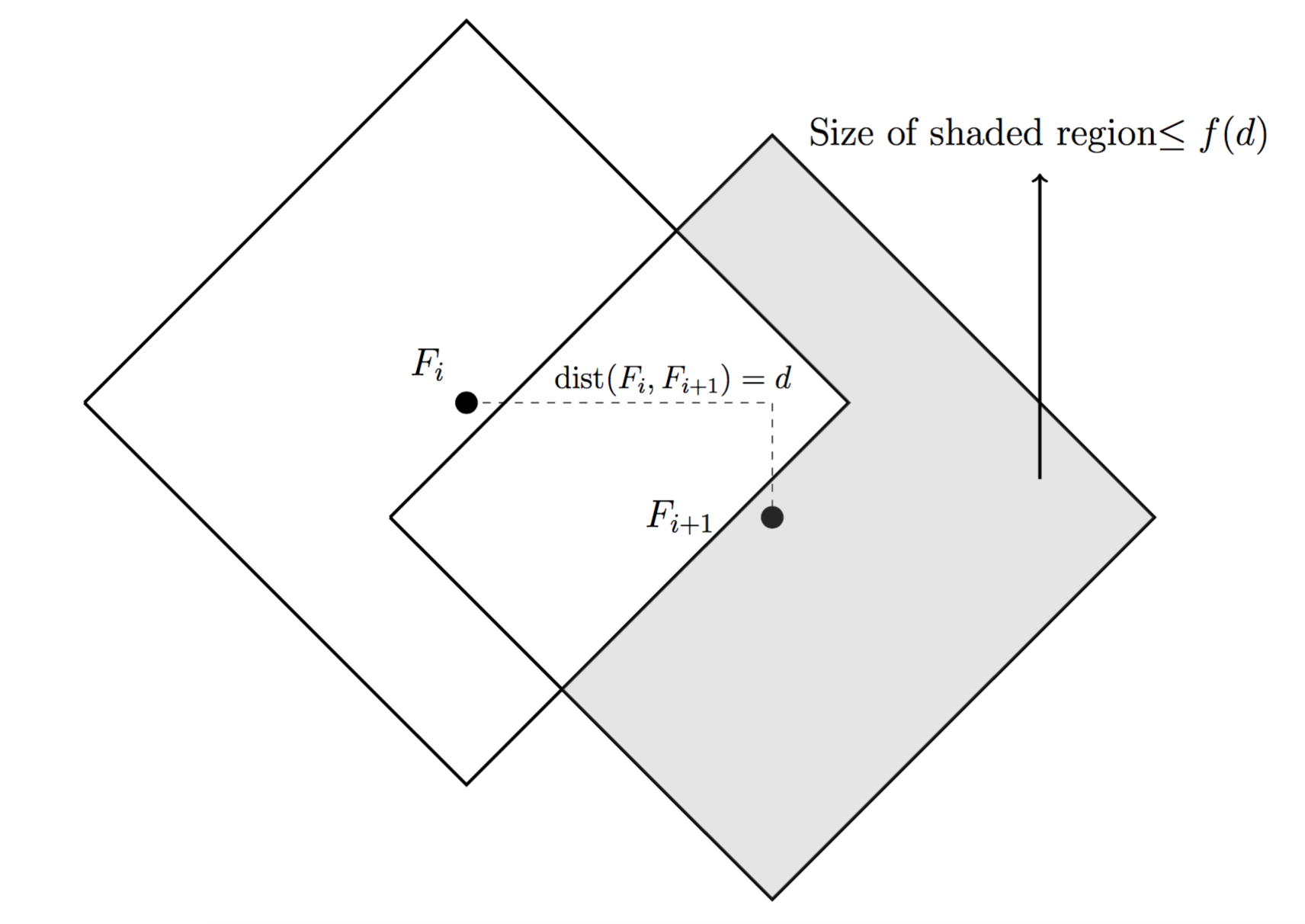}
	\end{center}
	\caption{Geometric interpretation of function $f(\cdot)$}
	\label{fig:overlap}
\end{wrapfigure}

We show in this section that one can not minimize path length $L$ and stop count $T$ simultaneously in the bi-objective CPPG. To minimize stop count (Corollary \ref{coro:1}), the overlap in coverage region for consecutive stops is minimized, resulting in longer path lengths. To minimize path length (Corollary \ref{coro:2}), traversals across the region are minimized which increases the overlap of coverage regions, resulting in more stops.

In the RC-CPPG and other CPPG variants, the parameter pair $(L,T)$ is feasible if there exists a covering path with at most $T$ stops and path length at most $L$. The set of all feasible pairs forms the feasible region. Note that if $(L,T)$ is a feasible pair, then both $(L,T+1)$ and $(L+\Delta L,T)$ are feasible pairs, where $\Delta L>0$. We explore the functional form of the trade-off between $L$ and $T$ to characterize the boundary of the feasible region. We define a trade-off function and use this function to solve the optimization version of the RC-CPPG. Our analysis of the boundary in the form of a trade-off between path length and stop count uses two functions based on $L$ and $T$:
\begin{itemize}
	\item the average distance between consecutive stops on a path of length $L$ with $T$ stops, $d=\frac{L}{T-1}$. Intuitively, when $d$ is large, the overlap between regions covered by consecutive stops is small which is associated with fewer stops.
	
	\item $f(d)=d(2k-\frac{d}{2})$ represents the maximum area of the region covered by a stop that is not covered by the previous stop on the path.
\end{itemize}

The function $f(\cdot)$ is an approximate measure of the area of the unique coverage region for consecutive stops. To minimize the number of stops $T$, one strives to maximize this area. Thus, this function plays an important role in the trade-off analysis governed by the choice of $L$ and $T$ through $d=\frac{L}{T-1}$.

\noindent We show the geometric interpretation of function $f(\cdot)$ in Figure \ref{fig:overlap}. The region covered by any stop is a diamond under the $l_1$ norm. If two consecutive stops $F_i$ and $F_{i+1}$ are separated by distance $d$, $f(d)$ is an upperbound of the area of the region covered by $F_{i+1}$ but not $F_i$ (shaded region in Figure \ref{fig:overlap}). Therefore, $f(d)$ can be interpreted as the maximum area covered by stop $F_{i+1}$ but not by the union of stops with lower index. If $d>2k$, then $F_i$ and $F_{i+1}$ serve disjoint regions and $f(d)=2k^{2}$, which is the size of the diamond region covered by $F_{i+1}$.

Theorem \ref{to-cont-R} characterizes the trade-off between $L$ and $T$.
\begin{theorem}[\textbf{Trade-off constraint for the RC-CPPG}\label{to-cont-R}]
	
	If $(L,T)$ is a feasible pair for the RC-CPPG with $T>1$, then
	\begin{equation}\label{eqn:trade-off-R}\tag{trade-off constraint}
		(T-1)f\Big(\frac{L}{T-1}\Big)\geq N-2k^{2},
	\end{equation}
	where $f(\cdot)$ is a function of the average distance between consecutive stops, $d$, defined as:
	\begin{equation}
		f(d)=
		\begin{cases}\label{def:f}
			d(2k-\frac{d}{2})& \text{if}~~d\in (0,2k]\\
			2k^{2}& \text{if}~~d\in (2k,\infty),
		\end{cases}
	\end{equation}
	and $N=mn$ is the area of rectangle $D$.
	
	Moreover, when $\frac{L}{T-1}\leq 2k$ (which is shown in Section \ref{sec:4.3} to hold for an optimal path), the \ref{eqn:trade-off-R} is equivalent to
	\begin{equation}\label{eqn:trade-off-R-ref}
		2kL-\frac{L^{2}}{2(T-1)}\geq N-2k^{2}.
	\end{equation}
\end{theorem}

\emph{Proof of Theorem \ref{to-cont-R}.}
Let $F_1-F_2-\cdots-F_T$ be a covering path where $\{F_i\}_{i=1}^{T}$ is the set of stops. Denote by $d_i$ the distance between $F_i$ and $F_{i+1}$. The total path length is $L=\sum_{i=1}^{T-1}d_i$. Let $S_i$ be the region covered by $F_i$ and $|S_i|$ its area. Note that $N\leq |\cup_{i=1}^{T}S_i|\leq |S_{1}|+\sum_{i=1}^{T-1}|S_{i+1}-S_i|$. We first show that $|S_{i+1}-S_i|\leq f(d_i)$ through the following lemma.
\begin{lemma}\label{coverlemma}
	For $k>0$ and $(a,b)\in \mathbb{R}^{2}$, let $\mathbb{B}\big((a,b),k\big)=\{(x,y)~|~|x-a|+|y-b|\leq k\}$ be the diamond region covered by $(a,b)\in \mathbb{R}^{2}$ with radius $k$. For $p, q\in \mathbb{R}$,
	\begin{equation}\label{eqn:cover}
		|\mathbb{B}\big((0,0),k\big)\cap \mathbb{B}\big((|p|+|q|,0),k\big)|\leq |\mathbb{B}\big((0,0),k\big)\cap \mathbb{B}\big((p,q),k\big)|.
	\end{equation}
\end{lemma}
\emph{Proof of Lemma \ref{coverlemma}}.
WLOG, we assume that $p\geq q\geq 0$ (else, we replace $(p,q)$ with $(|p|,|q|)$ and $(|p|,|q|)$ with $(|q|,|p|)$ if necessary, neither operations changes $|\mathbb{B}\big((0,0),k\big)\cap \mathbb{B}\big((p,q),k\big)|$). We now show that
\begin{equation}\label{eqn:subset}
	\mathbb{B}\big((0,0),k\big)\cap \mathbb{B}\big((|p|+|q|,0),k\big)\subseteq \mathbb{B}\big((0,0),k\big)\cap \mathbb{B}\big((p,q),k\big),
\end{equation}
and therefore prove inequality \eqref{eqn:cover}.

Since $p\geq q\geq 0$, $|p|+|q|=p+q$. For $(x_1,y_1)\in \mathbb{B}\big((0,0),k\big)\cap \mathbb{B}\big((p+q,0),k\big)$, we show that $(x_1,y_1)\in \mathbb{B}\big((p,q),k\big)$ and therefore $(x_1,y_1)\in \mathbb{B}\big((0,0),k\big)\cap \mathbb{B}\big((p,q),k\big)$.

If $x_1\leq p$, then $|x_1-p|+|y_1-q|\leq (p-x_1)+(q+|y_1|)\leq |x_1-p-q|+|y_1|\leq k$.

If $x_1>p$, then $|x_1-p|+|y_1-q|\leq (x_1-p)+(q+|y_1|)=(x_1+|y_1|)+(q-p)\leq |x_1|+|y_1|\leq k$.

Therefore, $|x_1-p|+|y_1-q|\leq k$ and $(x_1,y_1)\in \mathbb{B}\big((p,q),k\big)$. Thus \eqref{eqn:subset} follows, implying \eqref{eqn:cover}.
\qed

When $d_{i}$ is fixed, Lemma \ref{coverlemma} implies that the most efficient way to minimize $|S_{i}\cap S_{i+1}|$ is to locate $F_i$ and $F_{i+1}$ either vertically or horizontally within the grid.

From inequality \eqref{eqn:cover}, $|S_{i}\cap S_{i+1}|$ must be at least $|\mathbb{B}\big((0,0),k\big)\cap \mathbb{B}\big((d_i,0),k\big)|$. Together with the fact that $|S_{i+1}|=2k^{2}$, we have,
\begin{equation}\label{eqn:5}
	|S_{i+1}-S_i|=|S_{i+1}|-|S_{i}\cap S_{i+1}|\leq 2k^{2}-|\mathbb{B}\big((0,0),k\big)\cap \mathbb{B}\big((d_i,0),k\big)|.
\end{equation}
For $0\leq d_{i}<2k$, $\mathbb{B}\big((0,0),k\big)\cap \mathbb{B}\big((d_i,0),k\big)$ is a diamond region centered at $\big(\frac{d_i}{2},0\big)$ of radius $k-\frac{d_i}{2}$ with area $2\big(k-\frac{d_i}{2}\big)^{2}$. With $f(d_{i})$ defined in \eqref{def:f} we obtain
\[|\mathbb{B}\big((0,0),k\big)\cap \mathbb{B}\big((d_i,0),k\big)|=2\Big(k-\frac{d_i}{2}\Big)^{2}=2k^{2}-d_i\Big(2k-\frac{d_i}{2}\Big)=2k^{2}-f(d_i).\]
When $d_{i}\geq 2k$, $|\mathbb{B}\big((0,0),k\big)\cap \mathbb{B}\big((d_i,0),k\big)|=0$. Thus $f(d_i)=2k^{2}$.

Inequality \eqref{eqn:5} is then equivalent to $|S_{i+1}-S_i|\leq f(d_i)$. We thus obtain
\begin{equation}\label{eqn:second-to-last}
	N\leq |S_1|+\sum_{i=1}^{T-1}|S_{i+1}-S_{i}| \leq 2k^2+\sum_{i=1}^{T-1}f(d_i).
\end{equation}
Since $f(\cdot)$ is a concave function (see \eqref{def:f}), we have
\begin{equation}
	\begin{split}
		N-2k^{2}&\leq \sum_{i=1}^{T-1}f(d_i)\\
		&\leq (T-1)f\Big(\frac{\sum_{i=1}^{T-1}d_i}{T-1}\Big)~~\big(\textrm{from the concavity of}~f(\cdot)\big)\\
		&=(T-1)f\Big(\frac{L}{T-1}\Big).
	\end{split}
\end{equation}
Note that if $\frac{L}{T-1}\leq 2k$, $(T-1)f(\frac{L}{T-1})=(T-1)\cdot \frac{L}{T-1}(2k-\frac{L}{2(T-1)})=2kL-\frac{L^{2}}{2(T-1)}$. In this case, the \ref{eqn:trade-off-R} is equivalent to $2kL-\frac{L^{2}}{2(T-1)}\geq N-2k^{2}$. Thus \eqref{eqn:trade-off-R-ref} follows.
\qed

Theorem \ref{to-cont-R} provides a lower bound on the feasible region for $L$ and $T$ which we use to show the feasible paths defined in Section \ref{sec:4.2} are near-optimal. We show the trade-off constraint is almost tight in the sense that given a pair $(L,T)$ that satisfies the \ref{eqn:trade-off-R} at equality we can always find a feasible pair $(L^{'},T^{'})$ where $L$ is close to $L^{'}$ and $T$ is close to $T^{'}$.

\subsection{Near-optimal paths}\label{sec:4.2}
We define a group of paths called ``up-and-down paths''. One can consider the up-and-down path as a special case of the swath path (\cite{daganzo1984length}) which is shown to be near-optimal for the TSP in zones of different shapes. In defining an up-and-down path, we use the term ``traversal'' to represent the vertical line connecting points $(s,0)$ and $(s,m)$. An up-and-down path connects a set of traversals and the separation between consecutive traversals is a function of $d$ as in Definition \ref{def:UAD}. Once the path is defined, stops are located as described in Definition \ref{def:UAD}.

\begin{wrapfigure}{r}{0.5\textwidth}
	\begin{center}
		\includegraphics[width=0.48\textwidth]{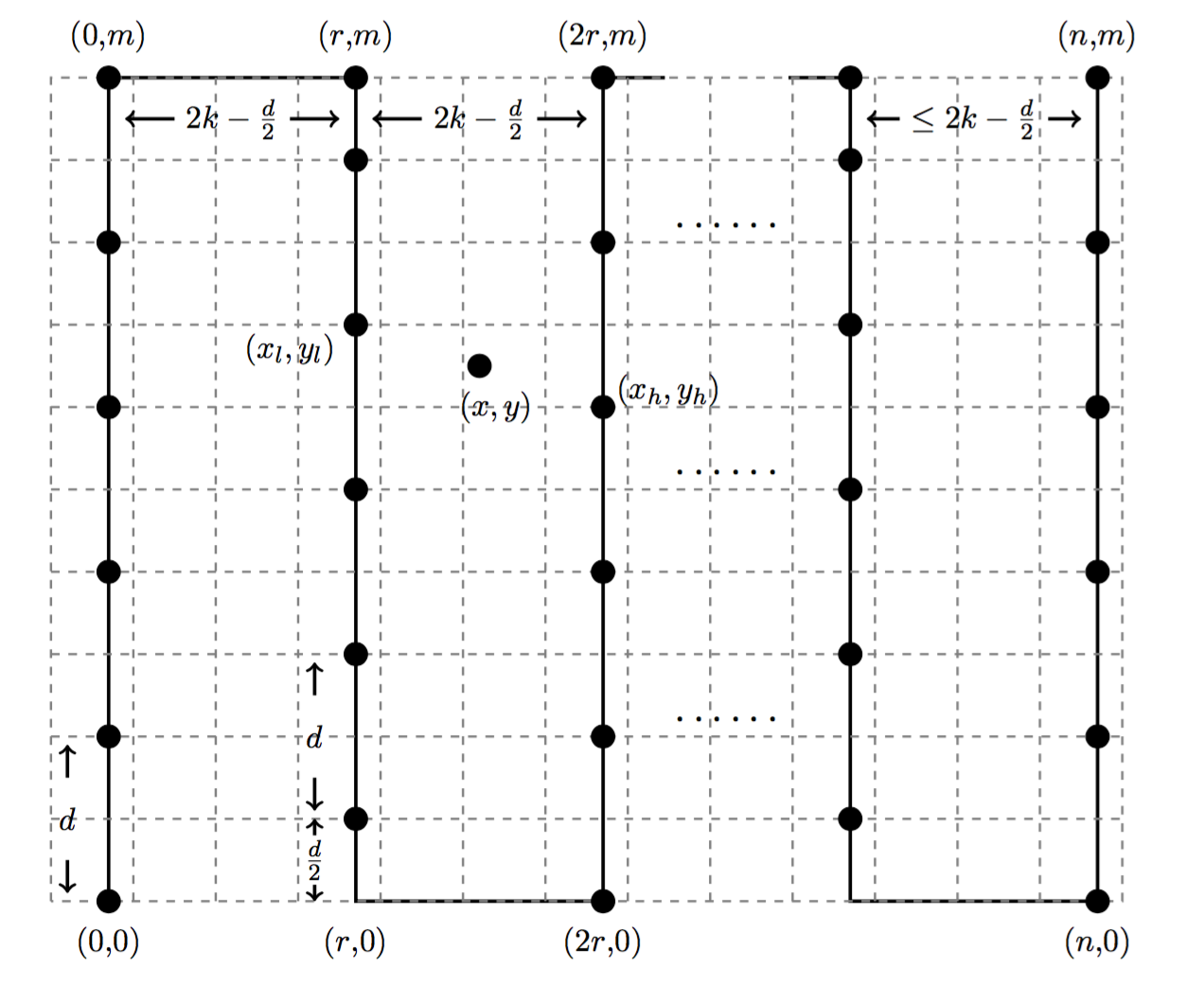}
	\end{center}
	\caption{Up-and-down path for the RC-CPPG}
	\label{fig:UAD-cont}
\end{wrapfigure}

\begin{definition}[\textbf{Type-$d$ up-and-down path}\label{def:UAD}]
	For $d\in (0,2k]$, define $r=2k-\frac{d}{2}$. In the RC-CPPG, a type-$d$ up-and-down path connects the following points sequentially (as shown in Figure \ref{fig:UAD-cont}):
	\[(0,0)\rightarrow(0,m)\rightarrow(r,m)\rightarrow(r,0)\rightarrow(2r,0)\rightarrow(2r,m)\rightarrow(3r,m)\rightarrow(3r,0)\rightarrow\cdots.\]
	
	For an odd traversal connecting $(2ir,0)$ to $(2ir,m)$ for $i=0,1,\cdots (i\leq \frac{n}{2r})$, we establish stops at $(2ir,jd)$ for $j=0,1,\cdots$ for $jd\leq m$. For an even traversal connecting $\big((2i+1)r,m\big)$ to $\big((2i+1)r,0\big)$ for $i=0,1,\cdots(i\leq \frac{n}{2r})$, we establish stops at $\big((2i+1)r,jd+\frac{d}{2}\big)$ for $j=0,1,\cdots$ for $jd+\frac{d}{2}\leq m$. Finally, we establish a stop at point $(\cdot,m)$ for each traversal to ensure coverage (which may not always be necessary as discussed in the proof of Proposition \ref{prop:feas-UAD}).
\end{definition}
Figure \ref{fig:UAD-cont} illustrates a type-$d$ up-and-down path where the black dots are the locations of selected stops. The points $(x,y), (x_h,y_h), (x_l,y_l)$ are used in the proof of Proposition \ref{prop:feas-UAD}.

Proposition \ref{prop:feas-UAD} shows the feasibility of the up-and-down path and Proposition \ref{prop:cost-UAD} computes the corresponding costs.
\begin{proposition}[\textbf{Feasibility of the up-and-down path}]\label{prop:feas-UAD}
	For any point $(x,y)\in D$, there exists a stop that covers $(x,y)$ on a type-$d$ up-and-down path.
\end{proposition}
\emph{Proof of Proposition \ref{prop:feas-UAD}.}
Assume that $(x,y)$ lies between traversals $i$ and $i+1$. Let $(x_h,y_h)$ be the highest stop on these two traversals with $y_h\leq y$ and $(x_l,y_l)$ be the lowest stop on these two traversals with $y_l\geq y$ (see Figure \ref{fig:UAD-cont} for illustration). From the alternating pattern of stop locations, we can pick $(x_h,y_h)$ and $(x_l,y_l)$ such that they are on adjacent traversals. Since the separation between traversals $i$ and $i+1$ is at most $2k-\frac{d}{2}$ (equal to $2k-\frac{d}{2}$ except for the rightmost one), we have $|x_h-x_l|\leq 2k-\frac{d}{2}$. With the alternating pattern of stop locations on traversals $i$ and $i+1$, we have $|y_h-y_l|\leq \frac{d}{2}$.

Since $y_h\leq y\leq y_l$ and $x$ is always between $x_h$ and $x_l$, we have
\[||(x,y)-(x_h,y_h)||_{1}+||(x,y)-(x_l,y_l)||_{1}=||(x_h,y_h)-(x_l,y_l)||_{1}\leq \Big(2k-\frac{d}{2}\Big)+\frac{d}{2}=2k.\]
This implies that $(x_h,y_h)$ or $(x_l,y_l)$ (and possibly both) covers $(x,y)$.
\qed

\begin{proposition}[\textbf{Cost of up-and-down path for the RC-CPPG}\label{prop:cost-UAD}]
	For $d\in (0,2k]$, the path length of a type-$d$ up-and-down path is at most $\frac{mn}{2k-d/2}+3m$ and the stop count is at most $\big(\frac{n}{2k-d/2}+2\big)\big(\frac{m}{d}+2\big)$. As $d\rightarrow 0$, the path length approaches $\frac{mn}{2k}+3m=\frac{mn}{2k}+O(m)$; if $d=2k$, the stop count is at most $(\frac{n}{k}+2)(\frac{m}{2k}+2)=\frac{mn}{2k^{2}}+O(m)$.
\end{proposition}
\emph{Proof of Propsition {\ref{prop:cost-UAD}}.}
In a type-$d$ up-and-down path, the separation between consecutive traversals (except the rightmost one) is $2k-\frac{d}{2}$, yielding at most $\frac{n}{2k-d/2}+2$ traversals. Two parts contribute to the path length: length from the traversals and length from traversal connections. The first part is bounded by $\big(\frac{n}{2k-d/2}+2\big)m$ and the second part is at most $n$. Since $m\geq n$, the total path length is at most $\big(\frac{n}{2k-d/2}+2\big)m+n\leq \frac{mn}{2k-d/2}+3m$.

The distance between consecutive stops on one traversal is $d$ (except for the topmost one).  Given that the number of stops on a single traversal is at most $\frac{m}{d}+2$, the total stop count is at most $\big(\frac{n}{2k-d/2}+2\big)\big(\frac{m}{d}+2\big)$.
\qed

\begin{theorem}[\textbf{Tightness result for the RC-CPPG}\label{tight-cont-R}]
	In the RC-CPPG, for any $(L,T)$ satisfying the \ref{eqn:trade-off-R} at equality, there exists a feasible up-and-down path of length $L^{'}$ with $T^{'}$ stops such that
	\begin{equation}\label{eqn:tight-L}
		L^{'}-L\leq 3m+2k,
	\end{equation}
	and
	\begin{equation}\label{eqn:tight-T}
		\frac{T^{'}}{T}\leq \Big(1+\frac{2k^{2}}{N-2k^2}\Big)\Big(1+\frac{4k}{n}\Big)^{2}.
	\end{equation}
\end{theorem}
Note that when $m,n$ are large, \eqref{eqn:tight-L} and \eqref{eqn:tight-T} suggest that both $\frac{L^{'}}{L}$ and $\frac{T^{'}}{T}$ can be arbitrarily close to 1. We show in Section \ref{sec:4.3} how this claim can be used
to obtain a near-optimal solution for the optimization problem.

\emph{Proof of Theorem \ref{tight-cont-R}.}
We discuss two cases based on the value of $d$: $d\leq 2k$ and $d>2k$. When $d\leq 2k$, we show that a type-$d$ up-and-down path is near-optimal; and when $d>2k$, a type-$2k$ up-and-down path is near-optimal.

\emph{Case 1: $d\leq 2k$.}

If $d=\frac{L}{T-1}\in (0,2k]$, the \ref{eqn:trade-off-R} at equality is $(T-1)f\big(\frac{L}{T-1}\big)=N-2k^{2}$. We can rewrite $L$ and $T$ as functions of $d$, \[T=T(d)=\frac{N-2k^{2}}{f(d)}+1=\frac{N-2k^2}{d(2k-d/2)}+1~~~(\textrm{from \eqref{def:f}})\]
and
\[L=L(d)=d(T(d)-1)=\frac{N-2k^2}{2k-d/2}.\]
Consider a type-$d$ up-and-down path and let $L^{'}$ and $T^{'}$ be the path length and stop count. From Proposition \ref{prop:cost-UAD} we have
\begin{equation}
	\begin{split}
		L^{'}\leq \frac{mn}{2k-d/2}+3m & =\frac{(N-2k^{2})+2k^{2}}{2k-d/2}+3m\\
		& \leq L(d)+3m+\frac{2k^{2}}{2k-d/2}\\
		& \leq L(d)+3m+2k.
	\end{split}
\end{equation}
For stop count, from Proposition \ref{prop:cost-UAD} we know
\begin{equation}
	\begin{split}
		T^{'} & \leq \Big(\frac{n}{2k-d/2}+2\Big)\Big(\frac{m}{d}+2\Big)\\
		& =\frac{n}{2k-d/2}\Big(1+\frac{2(2k-d/2)}{n}\Big)\cdot\frac{m}{d}\Big(1+\frac{2d}{m}\Big)\\
		& \leq\frac{N}{d(2k-d/2)}\Big(1+\frac{4k}{n}\Big)\Big(1+\frac{2d}{m}\Big)\\
		& \leq T(d)\Big(1+\frac{2k^{2}}{N-2k^2}\Big)\Big(1+\frac{4k}{n}\Big)\Big(1+\frac{4k}{n}\Big)~~~~~(\textrm{recall that}~d\leq 2k~\textrm{and}~m\geq n)\\
		& =T(d)\Big(1+\frac{2k^{2}}{N-2k^2}\Big)\Big(1+\frac{4k}{n}\Big)^{2}.
	\end{split}
\end{equation}

\emph{Case 2: $d>2k$.}

If $d=\frac{L}{T-1}>2k$, from the \ref{eqn:trade-off-R} we have
\[T\geq \frac{N-2k^{2}}{f(2k)}+1=T(2k), L\geq 2k(T-1)\geq L(2k). \]
From the analysis of \emph{Case 1}, the cost of a type-$2k$ up-and-down path satisfies \eqref{eqn:tight-L} and \eqref{eqn:tight-T} for $L=L(2k)$ and $T=T(2k)$. Since we are only increasing $L$ and $T$ in \emph{Case 2}, \eqref{eqn:tight-L} and \eqref{eqn:tight-T} still hold when $d>2k$.
\qed

\subsection{Optimization problem for the RC-CPPG}\label{sec:4.3}
The optimization version of the CPPG finds a feasible parameter pair $(L,T)$ that minimizes the objective function $C(L,T)$. This is equivalent to minimizing $C(L,T)$ over the feasible region defined by the \ref{eqn:trade-off-R}. Lemma \ref{convexity} guarantees the convexity of the feasible region, where the cost function $C(\cdot,\cdot)$ is increasing and convex.
\begin{lemma}\label{convexity}
	For given $N$ and $k$, the trade-off constraint $(T-1)f(\frac{L}{T-1})\geq N-2k^{2}$ defines a convex region.
\end{lemma}
\emph{Proof of Lemma \ref{convexity}.}
Since $f(\cdot)$ is concave, it is the minimum of a set of linear functions; i.e., $f(d)=\textrm{min}_i\{a_i d+b_i \}$. Therefore, the \ref{eqn:trade-off-R} is equivalent to $\textrm{min}_i\{ a_i L+b_i (T-1)\}\geq N-2k^{2}$. Note that $a_i L+b_i (T-1)\geq N-2k^{2}$ defines a halfspace in $\mathbb{R}^{2}$ for each $i$. Thus, $(T-1)f(\frac{L}{T-1})\geq N-2k^{2}$, which is the intersection of halfspaces, must be convex.
\qed

Given an increasing and convex function $C(L,T)$, consider the following subproblem:
\begin{equation}\label{eq:sub}
	\begin{split}
		\textrm{minimize}   &\quad C(L,T)\\
		\textrm{subject to} &\quad (T-1)f\Big(\frac{L}{T-1}\Big)\geq N-2k^{2}\\
		&\quad T\geq 2, L>0.
	\end{split}
\end{equation}
The problem minimizes a convex function over a convex set; therefore, it can be solved in polynomial time. Let $(L^{*},T^{*})$ be an optimal solution to the subproblem and $d^{*}=$min$\{ \frac{L^{*}}{T^{*}-1},2k\}$, with a corresponding type-$d^{*}$ up-and-down path. Theorem \ref{opt-R} provides a theoretical guarantee of the cost of up-and-down path for the case where $C(L,T)$ is linear.
\begin{theorem}\label{opt-R}
	If $C(L,T)=\alpha L+\beta T$ with $\alpha, \beta>0$ and $m\geq n\geq \frac{16k}{\varepsilon}$ for $\varepsilon\in (0,1)$, a type-$d^{*}$ up-and-down path provides a $(1+\varepsilon)$-approximation solution for the RC-CPPG.
\end{theorem}
\emph{Proof of Theorem \ref{opt-R}.}
Let $(L^{*},T^{*})$ be an optimal solution to subproblem \eqref{eq:sub}. Then $\alpha L^{*}+\beta T^{*}$ is a lower bound for the optimal function value. Since $(L^{*},T^{*})$ satisfies the \ref{eqn:trade-off-R} at equality, from Theorem \ref{tight-cont-R} we know there is a feasible pair $(L^{'},T^{'})$ such that $L^{'}-L^{*}\leq 3m+2k$ and $\frac{T^{'}}{T^{*}}\leq \big(1+\frac{2k^{2}}{N-2k^2}\big)\big(1+\frac{4k}{n}\big)^{2}$.

From \eqref{eqn:trade-off-R-ref}, we have $2kL^{*}-\frac{(L^{*})^{2}}{2(T^{*}-1)}\geq N-2k^{2}$ and $L^{*}\geq \frac{N-2k^{2}}{2k}=\frac{mn}{2k}-k$. Therefore,
\begin{equation}
	\begin{split}
		\frac{L^{'}}{L^{*}} & \leq 1+\frac{3m+2k}{L^{*}}\\
		& \leq 1+\Big(3m+2k\Big)\Big(\frac{2k}{mn-2k^{2}}\Big)\\
		& \leq 1+\frac{3m+2k}{8m/\varepsilon-k}~~~~~~~~~~\Big(\textrm{since}~~n\geq 16k/\varepsilon\Big)\\
		& \leq 1+\frac{5m}{7m/\varepsilon} \leq 1+\varepsilon.
	\end{split}
\end{equation}
From Theorem \ref{tight-cont-R}, we have
\begin{equation}
	\begin{split}
		\frac{T^{'}}{T^{*}} & \leq\Big(1+\frac{2k^{2}}{N-2k^2}\Big)\Big(1+\frac{4k}{n}\Big)^{2}\\
		& \leq\Big(1+\frac{2k^{2}}{256k^{2}/\varepsilon^{2}-2k^2}\Big)\Big(1+\frac{\varepsilon}{4}\Big)^{2}\\
		& \leq\Big(1+\frac{\varepsilon^{2}}{127}\Big)\Big(1+\frac{\varepsilon}{4}\Big)^{2}\leq 1+\varepsilon.
	\end{split}
\end{equation}
Thus $\alpha L^{'}+\beta T^{'}$ is at most $(1+\varepsilon)$ times the optimal solution because both $\frac{L^{'}}{L^{*}}$ and $\frac{T^{'}}{T^{*}}$ are at most $1+\varepsilon$.
\qed

With Theorems \ref{to-cont-R} and \ref{tight-cont-R}, we are able to solve two special cases where the objective only depends on one of the two costs; i.e., $C(L,T)=C(L)$ or $C(L,T)=C(T)$.
\begin{corollary}[Minimum path length in the RC-CPPG]\label{coro:1}
	If $C(L,T)=C(L)$, the optimal solution is achieved with $L^{*}=\frac{N}{2k}+O(m)$.
\end{corollary}
Let $(L,T)$ be a feasible pair. If $\frac{L}{T-1}\leq 2k$, inequality \eqref{eqn:trade-off-R-ref} implies $L\geq \frac{N-2k^{2}}{2k}=\frac{N}{2k}-k$. If $\frac{L}{T-1}>2k$, $L>2k(T-1)\geq 2k\cdot \frac{N-2k^{2}}{2k^{2}}=\frac{N}{k}-2k\geq \frac{N}{2k}+O(m)$. We can construct a type-$d$ up-and-down path that achieves this bound when $d\rightarrow 0$ (see Proposition \ref{prop:cost-UAD}).
\begin{corollary}[Minimum stop count in the RC-CPPG]\label{coro:2}
	If $C(L,T)=C(T)$, the optimal solution is achieved with $T^{*}=\frac{N}{2k^{2}}+O(m)$.
\end{corollary}
Since $f(d)\leq 2k^{2}$ for all $d>0$, we know from the \ref{eqn:trade-off-R} that
\begin{equation}
	T\geq \frac{N-2k^{2}}{f(\frac{L}{T-1})}+1\geq \frac{N-2k^{2}}{2k^{2}}+1=\frac{N}{2k^{2}}.
\end{equation}
We can construct a type-$2k$ up-and-down path that achieves the bound $\frac{N}{2k^{2}}+O(m)$ (see Proposition \ref{prop:cost-UAD}).

The optimal paths to minimize stop count and path length follow the same up-and-down pattern with different parameters. Note that we can not minimize path length and stop count simultaneously: the path length minimizing path aims to minimize traversals which requires more stops while the stop count minimizing path forms a tessellation which decreases the separation between traversals ($r=2k-\frac{d}{2}$) thus increasing path length. The structure of the optimal paths also coincides with that of the optimal solutions in \cite{carlsson2014continuous} when the access cost is measured by the $l_1$ norm.

In summary, we introduce a trade-off constraint to quantify the trade-off between path length $L$ and stop count $T$ in a covering path for the RC-CPPG. We construct a family of feasible paths that traverse the rectangle in an up-and-down pattern. For the optimization problem, we show the costs of the up-and-down paths match the lower bound derived from the trade-off constraint. The optimal solution can be found quickly through a simple algorithm based on a convex relaxation of the optimization problem. This simple approach used to solve the RC-CPPG is extended in Sections \ref{continuous} and \ref{discrete} to solve the C-CPPG and D-CPPG, respectively. With Proposition \ref{prop:1}, this then gives a complete solution approach to the CPPG.

\section{Analysis of the C-CPPG \label{continuous}}
In this section, we focus on the C-CPPG where the coverage radius $k$ is an integer and the rectangle $D$ is covered by stops selected from $D_{int}$. The C-CPPG is a special case of the RC-CPPG where the stop locations are chosen from $D_{int}$ rather than $D$. Therefore, the feasible region in the C-CPPG is a subset of that in the RC-CPPG and the \ref{eqn:trade-off-R} for the RC-CPPG is valid for the C-CPPG. We strengthen the \ref{eqn:trade-off-R} to provide a tighter lower bound for the feasible region in the C-CPPG. We then generalize the up-and-down path by mixing two types of up-and-down paths. We show the costs of the generalized paths match the tighter lower bound within a constant factor. Finally, we present a constant-factor approximation algorithm for solving the optimization problem with a linear objective function. Since the results follow directly from those in Section \ref{relax}, we present proofs in the Appendix and highlight the differences here.
\subsection{Stronger trade-off constraint for the C-CPPG}\label{sec:5.1}
Given that stop locations are restricted to integer points, the distance between consecutive stops $d$ must be integer. Hence, we strengthen $f(\cdot)$ with a piecewise-linear function $f_{LB-C}(\cdot)$ which connects all points on $f(\cdot)$ with integer input values; i.e.,

\begin{equation}\label{eqn:fLB-C}
	f_{LB-C}(d)=
	\begin{cases}
		(t+1-d)f(t)+(d-t)f(t+1)& \textrm{if}~~d\in [t,t+1), t\in [2k-1];\\
		2k^{2}& \textrm{if}~~d\geq 2k.
	\end{cases}
\end{equation}
Essentially, $f_{LB-C}(\cdot)$ is the piecewise-linear function connecting
\begin{equation}\label{eqn:boudary-LBC}
	\big(1,f(1)\big)\rightarrow\big(2,f(2)\big)\rightarrow\cdots\rightarrow\big(2k,f(2k)\big)\rightarrow\big(\infty,f(2k)\big).
\end{equation}
Since $f_{LB-C}(d)=f(d)$ when $d$ is integer, $f_{LB-C}(\cdot)$ also represents the maximum area of the region covered by a stop that is not covered by previous stops. From \eqref{def:f}, $f(\cdot)$ is a concave function. Therefore, $f_{LB-C}(\cdot)$ must be a concave function satisfying $f_{LB-C}(\cdot)\leq f(\cdot)$. This implies that constraint \eqref{eqn:trade-off-cont} is stronger than the \ref{eqn:trade-off-R} for the RC-CPPG.
\begin{theorem}[\textbf{Trade-off constraint for the C-CPPG}\label{to-cont-C}]
	In the C-CPPG, if $(L,T)$ is a feasible pair with $T>1$, then
	\begin{equation}\label{eqn:trade-off-cont}
		(T-1)f_{LB-C}\Big(\frac{L}{T-1}\Big)\geq N-2k^{2}.
	\end{equation}
	
	Moreover, the boundary of \eqref{eqn:trade-off-cont}, $(T-1)f_{LB-C}\big(\frac{L}{T-1}\big)=N-2k^{2}$ is a polyline connecting
	\begin{equation}\label{eqn:polyline-LBC}
		\begin{split}
			&\Big(\frac{N^{*}}{f(1)},\frac{N^{*}}{f(1)}\Big)\rightarrow\Big(\frac{2N^{*}}{f(2)},\frac{N^{*}}{f(2)}\Big)\rightarrow\cdots\rightarrow\Big(\frac{iN^{*}}{f(i)},\frac{N^{*}}{f(i)}\Big)\rightarrow\\
			&\cdots\rightarrow\Big(\frac{2kN^{*}}{f(2k)},\frac{N^{*}}{f(2k)}\Big)\rightarrow\Big(\infty,\frac{N^{*}}{f(2k)}\Big),
		\end{split}
	\end{equation}
	in the $(L,T)$ plane, where $N^{*}=N-2k^{2}$.
\end{theorem}

Compared to Theorem \ref{to-cont-R}, we replace the function $f(\cdot)$ with $f_{LB-C}(\cdot)$ to strengthen the trade-off constraint. The boundary of the feasible region in the C-CPPG, defined by \eqref{eqn:polyline-LBC}, connects several points on \eqref{eqn:trade-off-R-ref}. For both the trade-off constraint and the boundary of the feasible region, we replace a smooth function with a piecewise-linear one. Lemma \ref{lem:piecewise} establishes a basic property of the piecewise-linear function which is useful in the study of the boundary of \eqref{eqn:trade-off-cont}.
\begin{lemma}\label{lem:piecewise}
	Let $0<a_1<\cdots<a_{n}<a_{n+1}=\infty$ and $0<b_1<\cdots<b_{n}=b_{n+1}$ be nonnegative increasing sequences. Let $g(\cdot)$ be a piecewise-linear function defined on $[a_1,\infty)$ corresponding to the polyline connecting $(a_1,b_1)\rightarrow(a_2,b_2)\rightarrow\cdots\rightarrow(a_n,b_n)\rightarrow(a_{n+1},b_{n+1})$. If $g(\cdot)$ is a concave function, for any constant $C>0$, $Y\cdot g(\frac{X}{Y})=C$ is a polyline connecting the following points
	\begin{equation}\label{eqn:polyline-C}
		\Big(\frac{a_{1}C}{b_1},\frac{C}{b_1}\Big)\rightarrow\Big(\frac{a_{2}C}{b_2},\frac{C}{b_2}\Big)\rightarrow\cdots\rightarrow\Big(\frac{a_{n}C}{b_n},\frac{C}{b_n}\Big)\rightarrow\Big(\frac{a_{n+1}C}{b_{n+1}},\frac{C}{b_{n+1}}\Big).
	\end{equation}
	Moreover, $Yg(\frac{X}{Y})$ is a convex function and \eqref{eqn:polyline-C} corresponds to a convex piecewise-linear function.
\end{lemma}	

For the boundary of \eqref{eqn:trade-off-cont}, note that $f_{LB-C}(\cdot)$ is a piecewise-linear concave function (derived from the concavity of $f(\cdot)$). From Lemma \ref{lem:piecewise} and \eqref{eqn:boudary-LBC} we know $(T-1)f_{LB-C}(\frac{L}{T-1})=N-2k^{2}$ is equivalent to polyline \eqref{eqn:polyline-LBC}.
\subsection{Near-optimal mixed up-and-down path for the C-CPPG}
When stops are restricted to integer points, the up-and-down path defined in Section \ref{relax} is not sufficient because the distances between consecutive stops, $d$, and the separation between traversals, $2k-\frac{d}{2}$, may not be integers. Having obtained the stronger trade-off constraint for the C-CPPG in Theorem \ref{to-cont-C}, we now define mixed up-and-down paths that are then shown to be close to any feasible parameter pair on the boundary of \eqref{eqn:trade-off-cont}. The mixed up-and-down path is conducted with a part of the rectangle covered by a type-$d$ up-and-down path as in Section \ref{relax} and a part of the rectangle covered by a type-$(d+2)$ up-and-down path. We restrict $d$ to be an even integer in $[2k]$ for each path that comprises the mixed path scheme.

\begin{definition}[\textbf{Mixed up-and-down path}]\label{def:MUAD}
	In the C-CPPG, for even values of $d\in [2k-2]$ and $\gamma\in [0,1)$, we divide the $m\times n$ rectangle into two rectangles of sizes $\left\lceil \gamma n\right\rceil \times m$ and $\left\lceil (1-\gamma)n\right\rceil \times m$. A type-$(d,\gamma)$ mixed up-and-down path covers the $\left\lceil \gamma n\right\rceil \times m$ rectangle with a type-$d$ up-and-down path and the $\left\lceil (1-\gamma)n\right\rceil \times m$ rectangle with a type-$(d+2)$ up-and-down path, and connects the two paths at the common boundary of the two rectangles.
\end{definition}

For even values of $d$, both $d$ and $2k-\frac{d}{2}$ are integers and all stops in the mixed up-and-down path are located at integer points. The mixed up-and-down path is always feasible since both parts of the rectangle are covered by the type-$d$ and type-$(d+2)$ up-and-down paths, respectively. Proposition \ref{prop:cost-MUAD} computes the costs of a mixed up-and-down path. In Section \ref{sec:5.3}, we discuss the selection of $\gamma$.
\begin{proposition}[\textbf{Cost of mixed up-and-down path}\label{prop:cost-MUAD}]
	For any even $d\in [2k-2]$ and $\gamma\in [0,1)$, the path length $L$ of a type-$(d,\gamma)$ up-and-down path is at most $\frac{\gamma mn}{2k-d/2}+\frac{(1-\gamma) mn}{2k-(d+2)/2}+10m$ and the stop count $T$ is at most $\frac{\gamma mn}{d(2k-d/2)}+\frac{(1-\gamma) mn}{(d+2)(2k-(d+2)/2)}+10m+12$.
\end{proposition}
The mixed up-and-down path can be divided into three parts: the type-$d$ up-and-down path, the type-$(d+2)$ up-and-down path and the segment connecting these two paths. We compute the costs of each part based on Proposition \ref{prop:cost-UAD}.

Based on Proposition \ref{prop:cost-MUAD}, Theorem \ref{tight-cont-C} provides an approximate upper bound for the feasible region in the C-CPPG.
\begin{theorem}[\textbf{Tightness result for the C-CPPG}\label{tight-cont-C}]
	Let $f_{UB-C}(\cdot)$ be the piecewise-linear function defined on $[2,\infty)$ that corresponds to the polyline connecting
	\[\big(2,f(2)\big)\rightarrow\big(4,f(4)\big)\rightarrow\cdots\rightarrow\big(2i,f(2i)\big)\rightarrow\big(2k,f(2k)\big)\rightarrow\big(\infty,f(2k)\big).\]
	For any $(L,T)$ satisfying
	\begin{equation}\label{eqn:trade-off-C-eq}
		(T-1)f_{UB-C}\Big(\frac{L}{T-1}\Big)=N-2k^{2},
	\end{equation}
	there exists a feasible pair $(L^{'},T^{'})$ derived from a mixed up-and-down path such that
	\begin{equation}\label{eqn:23}
		L^{'}-L\leq \frac{2k^{2}}{N-2k^{2}}L+10m,
	\end{equation}
	and
	\begin{equation}\label{eqn:24}
		T^{'}-{T}\leq \frac{2k^{2}}{N-2k^{2}}T+10m+12.
	\end{equation}
	Moreover, \eqref{eqn:trade-off-C-eq} is a polyline connecting
	\begin{equation}\label{eqn:poly-UBC}
		\begin{split}
			&\Big(\frac{2N^{*}}{f(2)},\frac{N^{*}}{f(2)}\Big)\rightarrow\Big(\frac{4N^{*}}{f(4)},\frac{N^{*}}{f(4)}\Big)\rightarrow\cdots\rightarrow\Big(\frac{2iN^{*}}{f(2i)},\frac{N^{*}}{f(2i)}\Big)\rightarrow\\
			&\cdots\rightarrow\Big(\frac{2kN^{*}}{f(2k)},\frac{N^{*}}{f(2k)}\Big)\rightarrow\Big(\infty,\frac{N^{*}}{f(\infty)}\Big),
		\end{split}
	\end{equation}
\end{theorem}

When $m$ and $n$ are sufficiently large, $\frac{2k^{2}}{N-2k^{2}}$ can be arbitrarily close to 0. Therefore, \eqref{eqn:poly-UBC} serves as an approximate upper bound of the feasible region in the C-CPPG.
\subsection{Optimization problem for the C-CPPG}\label{sec:5.3}
Given a convex and increasing function $C(L,T)$, we solve the following convex subproblem:
\begin{equation}\label{eq:sub-C}
	\begin{split}
		\textrm{minimize}   &\quad C(L,T)\\
		\textrm{subject to} &\quad (T-1)f_{LB-C}\Big(\frac{L}{T-1}\Big)=N-2k^{2}\\
		&\quad T\geq 2, L>0.
	\end{split}
\end{equation}
Let $(L^{*}, T^{*})$ be an optimal solution to \eqref{eq:sub-C} and $d^{*}=\frac{L^{*}}{T^{*}}$. We first consider the case where $2\leq d^{*}< 2k$. Denote $\mathbb{S}=\{2,4,\cdots,2i,\cdots,2k\}$ to be the set of even numbers in $[2k]$. Let $l$ be the largest number in $\mathbb{S}$ that is no greater than $d^{*}$; let $r$ be the smallest number in $\mathbb{S}$ that is greater than $d^{*}$. Select $\gamma\in [0,1)$ such that $d^{*}=\gamma l+(1-\gamma)r$ (for example, if $d^{*}=3$ we have $l=2, r=4$ and $\gamma=\frac{1}{2}$). We use a type-$(l,\gamma)$ mixed up-and-down path to cover the rectangle. If $d^{*}<2$, we consider a type-2 up-and-down path; and if $d^{*}\geq 2k$ we consider a type-$2k$ up-and-down path. Note that while $d^{*}=\frac{L^{*}}{T^{*}}$ may not be an even integer, $d$ and $d+2$ used to create the mixed up-and-down path are always even integers.

\begin{wrapfigure}{r}{0.5\textwidth}
	\begin{center}
		\includegraphics[width=0.48\textwidth]{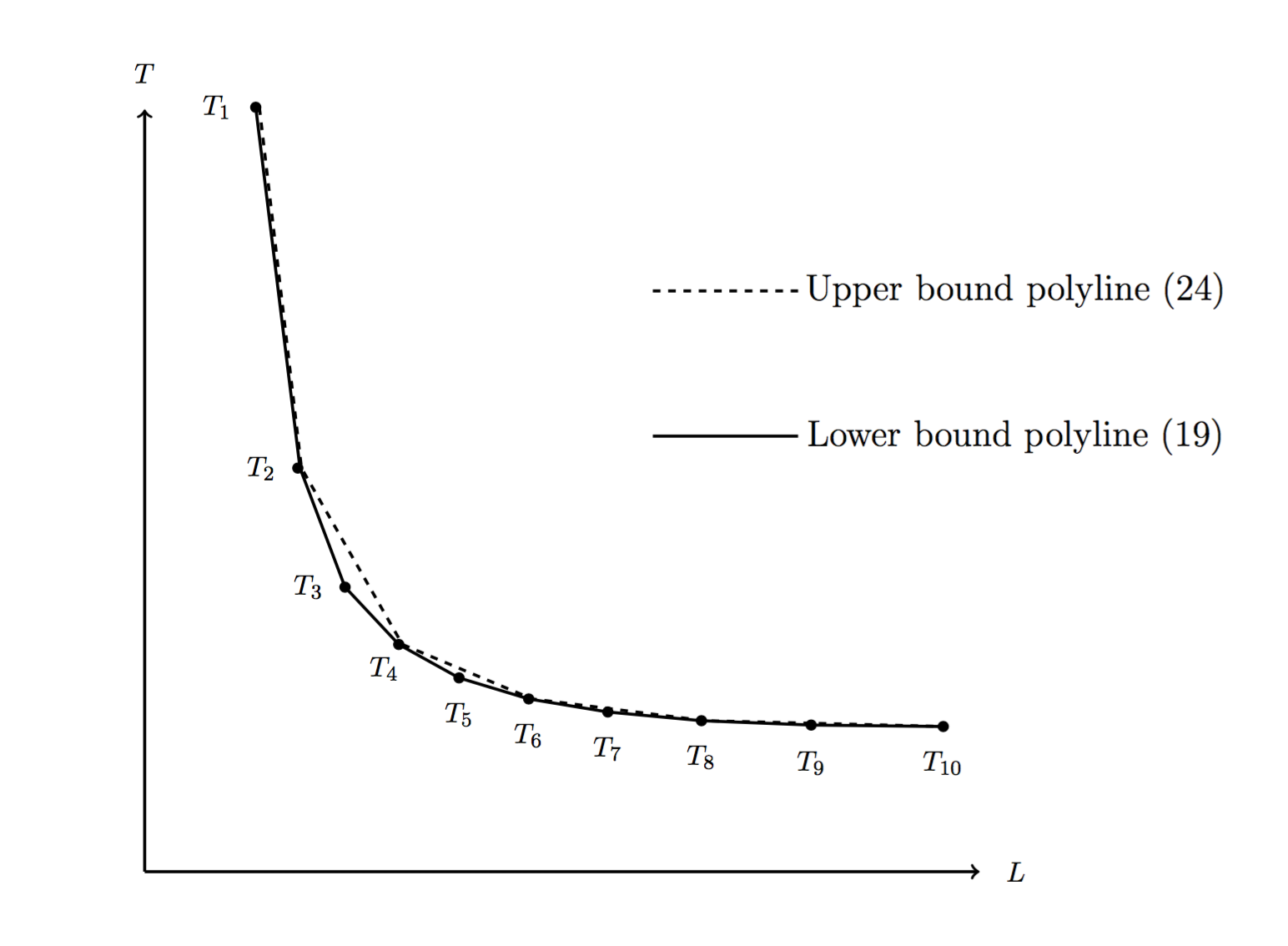}
	\end{center}
	\caption{Gap between lower bound and upper bound\\for the C-CPPG.}
	\label{fig:GAP}
\end{wrapfigure}

For the C-CPPG, \eqref{eqn:polyline-LBC} and \eqref{eqn:poly-UBC} provide a lower bound and an approximate upper bound, respectively, for the feasible region. Unlike the RC-CPPG, the bounds provided by \eqref{eqn:polyline-LBC} and \eqref{eqn:poly-UBC} are not tight. A gap exists because the turning points of \eqref{eqn:polyline-LBC} contain that of \eqref{eqn:poly-UBC} while there are turning points of \eqref{eqn:polyline-LBC} not contained in those of \eqref{eqn:poly-UBC}. Figure \ref{fig:GAP} illustrates the gap for $k=5$. The solid line connecting $T_1, T_2,\cdots, T_{10}$ is the lower bound polyline \eqref{eqn:polyline-LBC}; the dashed line connecting even turning points $T_2, T_4,\cdots, T_{10}$ represents the upper bound \eqref{eqn:poly-UBC}. Lemma \ref{approx-C} shows the relative gap between \eqref{eqn:polyline-LBC} and \eqref{eqn:poly-UBC} is at most $\frac{9}{8}$ for $k\geq 3$ ($\frac{3}{2}$ for $k=1$ and $\frac{7}{6}$ for $k=2$).

\begin{lemma}\label{approx-C}
	For any point $(L_1, T_1)$ on the lower bound polyline \eqref{eqn:polyline-LBC}, there exists point $(L_2, T_2)$ on the upper bound polyline \eqref{eqn:poly-UBC} such that $\frac{L_2}{L_1}\leq \frac{11}{10}$ for $k\geq 3$ ($\frac{3}{2}$ for $k=1$ and $\frac{7}{6}$ for $k=2$), $\frac{T_2}{T_1}\leq \frac{9}{8}$.
\end{lemma}

For the optimization problem with a linear objective function, we obtain a lower bound by solving subproblem \eqref{eq:sub-C}; i.e., minimizing the objective on the lower bound polyline \eqref{eqn:poly-UBC}. Since the polyline \eqref{eqn:poly-UBC} characterizes the costs of the mixed up-and-down paths, minimizing $C(L,T)$ on \eqref{eqn:poly-UBC} provides an approximate upper bound for the optimization problem. Lemma \ref{approx-C} implies that the lower bound and the upper bound for the optimizarion problem have a gap of at most $\frac{3}{2}$. This gap decreases as $k$ increases. Theorem \ref{opt-C} formally proves approximation results for minimizing a linear function $C(L,T)$.
\begin{theorem}\label{opt-C}
	If $C(L,T)=\alpha L+\beta T$ with $\alpha, \beta\geq 0$ and $m\geq n\geq \frac{100k^{2}}{\varepsilon}\geq \frac{100k}{\varepsilon}$ where $\varepsilon\in (0,1)$, the mixed up-and-down path provides a $(\frac{9}{8}+\varepsilon)$-approximation solution for $k\geq 3$ ($\frac{3}{2}+\varepsilon$ for $k=1$ and $\frac{7}{6}+\varepsilon$ for $k=2$).
\end{theorem}

\begin{corollary}[Minimum stop count in the C-CPPG]
	If $C(L,T)=C(T)$, the optimal solution is achieved by $T^{*}=\frac{N}{2k^{2}}+O(m)$.
\end{corollary}
\emph{Proof.}
The result here is the same as in Corollary \ref{coro:2} because the stop count minimizing path in the RC-CPPG, a type-$2k$ up-and-down path, is a covering path in the C-CPPG.
\qed

\textbf{Remark:} For the case where $C(L,T)=C(L)$, we only know that the optimal solution is between $\frac{N}{2k-1/2}+O(m)$ and $\frac{N}{2k-1}+O(m)$. A type-$1$ up-and-down path (which corresponds to the lower bound $\frac{N}{2k-1/2}+O(m)$ in this setting) is infeasible since not all stops are located at integer points.

\section{Analysis of the D-CPPG \label{discrete}}
In this section we focus on the D-CPPG where the coverage radius $k$ is an integer and both coverage region and stop locations are restricted to integer points in $D_{int}$. Detailed proofs are provided in the Appendix.
\subsection{Trade-off constraint for the D-CPPG}
As in Section \ref{sec:5.1}, we provide a piecewise-linear function $f_{LB-D}(\cdot)$ that represents the maximum number of integer points covered by a stop but not covered by previous stops. We define the function value of $f_{LB-D}(\cdot)$ with integer inputs:
\begin{equation}\label{eqn:trade-off-disc}
	f_{LB-D}(d)=
	\begin{cases}
		d(2k+1-\frac{d}{2})+\frac{1}{2}& \textrm{if}~~d\in[2k]~\textrm{is odd};\\
		d(2k+1-\frac{d}{2})& \textrm{if}~~d\in[2k]~\textrm{is even};\\
		2k^{2}+2k+1& \textrm{if}~~d=2k+1.
	\end{cases}
\end{equation}
For $d\in (t,t+1)$, where $t\in [2k]$, $f_{LB-D}(d)=(t+1-d)f_{LB-D}(t)+(d-t)f_{LB-D}(t+1)$.

Equivalently, $f_{LB-D}(\cdot)$ is the piecewise-linear function connecting
\begin{equation}\label{eqn:boudary-LBD}
	\big(1,f_{LB-D}(1)\big)\rightarrow\big(2,f_{LB-D}(2)\big)\rightarrow\cdots\rightarrow\big(2k+1,f_{LB-D}(2k+1)\big)\rightarrow\big(\infty,f_{LB-D}(2k+1)\big).
\end{equation}
Similar to Theorems \ref{to-cont-R} and \ref{to-cont-C}, we provide a trade-off constraint for the D-CPPG using function $f_{LB-D}(\cdot)$. There are two differences in the trade-off constraint for the D-CPPG, \eqref{eqn:to-disc}.

(1). We use a different piecewise-linear function $f_{LB-D}(\cdot)$. The function $f_{LB-D}(d)$ counts the maximum number of integer points covered by stop $F_{i+1}$ but not $F_{i}$ given that the two stops are at distance $d$.

(2). On the right hand side of the trade-off constraint \eqref{eqn:to-disc} we use $2k^2+2k+1$ instead of $2k^{2}$ in \eqref{eqn:trade-off-cont} since the number of integer points covered by one stop is $2k^2+2k+1$ in the D-CPPG.

\begin{theorem}[\textbf{Trade-off constraint for the D-CPPG}\label{tradeoff-D}]
	In the D-CPPG, if $(L,T)$ is a feasible pair with $T>1$, then
	\begin{equation}\label{eqn:to-disc}
		(T-1)f_{LB-D}\Big(\frac{L}{T-1}\Big)\geq N-(2k^2+2k+1).
	\end{equation}
	
	Moreover, the boundary of \eqref{eqn:to-disc}, $(T-1)f_{LB-D}(\frac{L}{T-1})=N-(2k^2+2k+1)$, is a piecewise-linear convex function connecting
	\begin{equation}\label{eqn:polyline-D}
		\begin{split}
			&\Big(\frac{N^{**}}{f_{LB-D}(1)},\frac{N^{**}}{f_{LB-D}(1)}\Big)\rightarrow\Big(\frac{2N^{**}}{f_{LB-D}(2)},\frac{N^{**}}{f_{LB-D}(2)}\Big)\rightarrow\cdots\\
			&\rightarrow\Big(\frac{(2k+1)N^{**}}{f_{LB-D}(2k+1)},\frac{N^{**}}{f_{LB-D}(2k+1)}\Big)\rightarrow\Big(\infty,\frac{N^{**}}{f_{LB-D}(2k+1)}\Big),
		\end{split}
	\end{equation}
	where $N^{**}=N-(2k^{2}+2k+1)$.
\end{theorem}

The fact that \eqref{eqn:polyline-D} is the boundary of \eqref{eqn:to-disc} follows from Lemma \ref{lem:piecewise} and the fact that $f_{LB-D}(\cdot)$ is a concave piecewise-linear function. The proof of Theorem \ref{tradeoff-D} follows the same structure as Theorem \ref{to-cont-C}, see Appendix \ref{sec:B1} for details.
\subsection{Near-optimal mixed discrete up-and-down path for the D-CPPG}
\begin{figure}
	\begin{center}
		\includegraphics[scale=0.33]{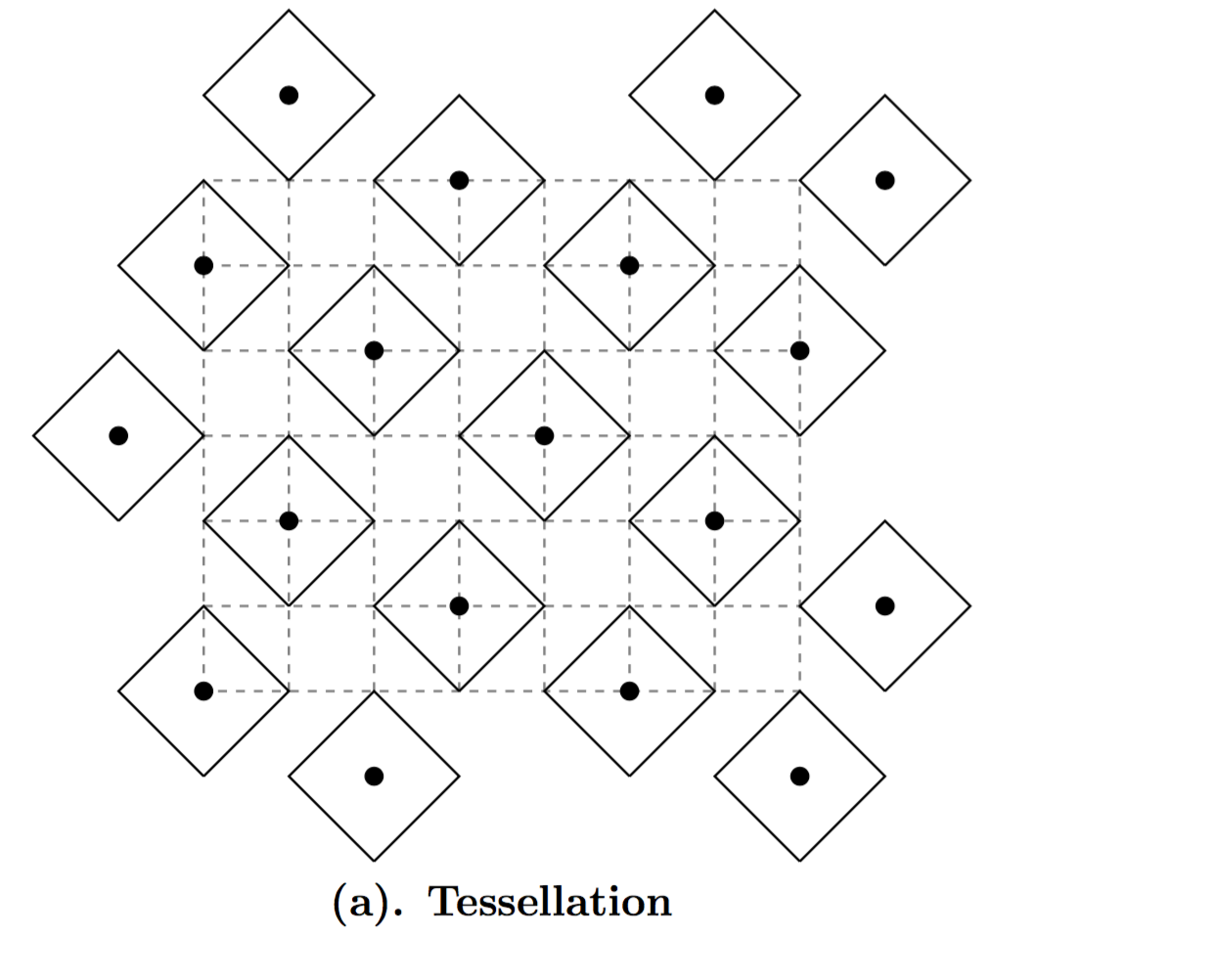}
		\includegraphics[scale=0.33]{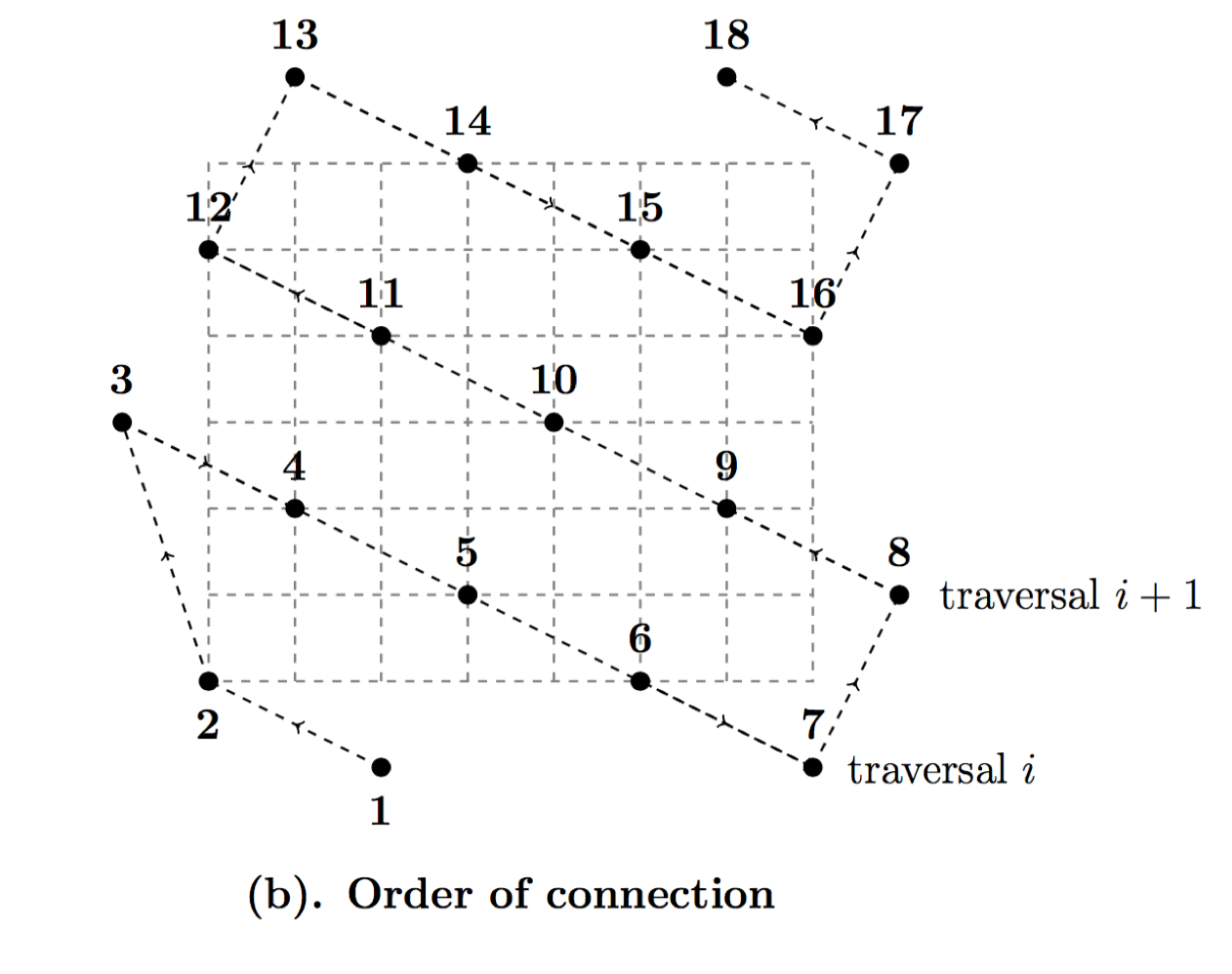}
		\includegraphics[scale=0.33]{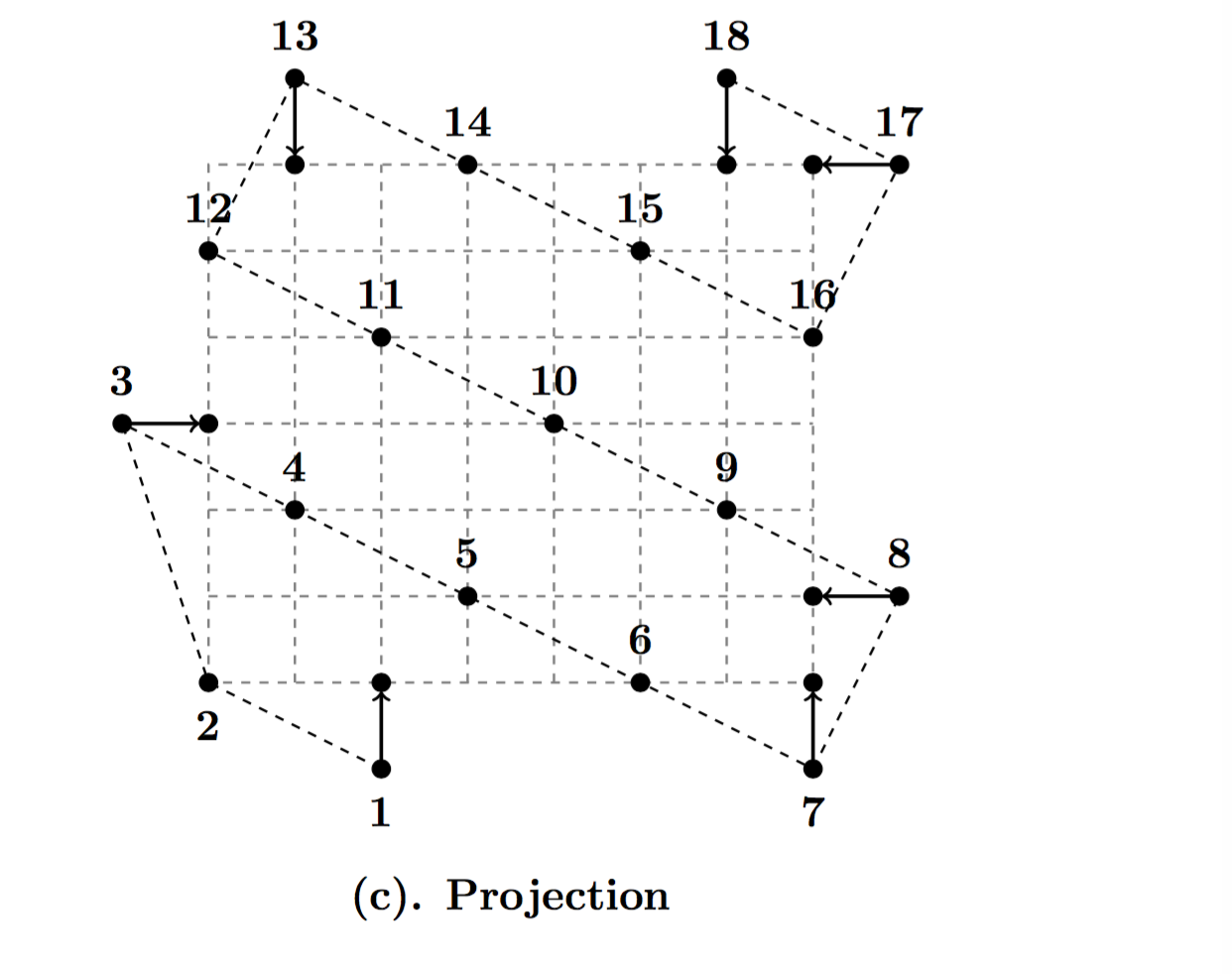}
		\includegraphics[scale=0.33]{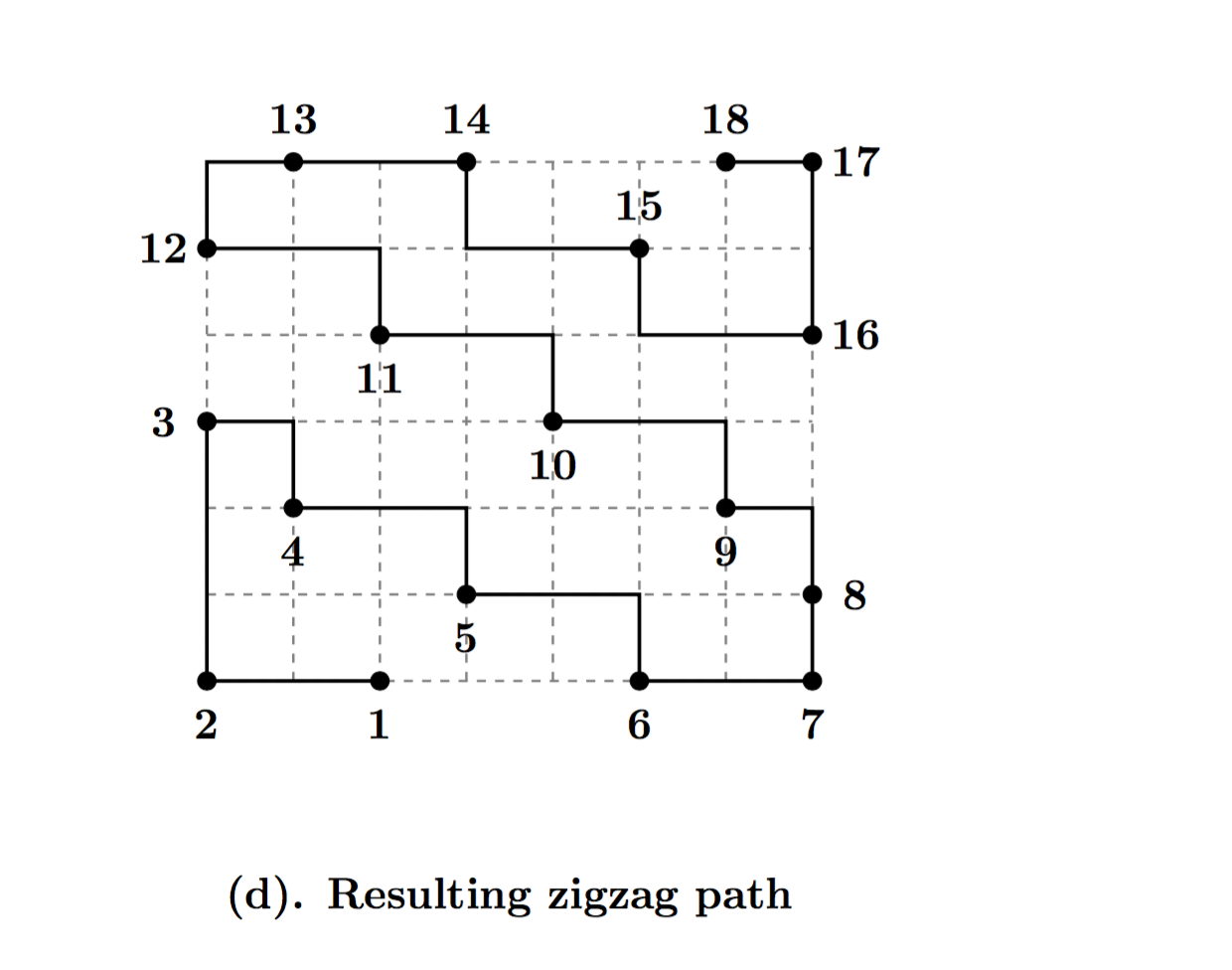}
		\caption{Stop location and path pattern of type-(2$k$+1) discrete up-and-down path (zigzag path)} \label{fig:DUAD}
	\end{center}
\end{figure}

Similar to Definition \ref{def:UAD}, we define the type-$d$ discrete up-and-down path for the D-CPPG. In Definition \ref{def:UAD}, the distance between consecutive stops along a traversal is $d$ and the separation between traversals is $2k-\frac{d}{2}$. As with the C-CPPG, we restrict $d$ to even integers, except $d=1$. For $d=1,2,4,6,\cdots,2k-2,2k$, the type-$d$ discrete up-and-down path traverses the coverage region $D_{int}$ in a up-and-down fashion but the separation between traversals changes as described in Definition \ref{def:DUAD}. For $d=2k+1$, we introduce the zigzag path which has a different pattern than the up-and-down path. We show that the zigzag path minimizes stop count in Corollary \ref{minstop} and that a suitable combination of zigzag path and up-and-down path provides a high quality approximation solution in Theorem \ref{approx-D}.
\begin{definition}[\textbf{Type-$d$ discrete up-and-down path}]\label{def:DUAD}
	We define the type-$d$ discrete up-and-down path in the D-CPPG for $d=1,2,4,6,\cdots,2k-2,2k,2k+1$.
	\begin{itemize}
		\item In a type-$1$ discrete up-and-down path, the distance between consecutive stops is 1 along each traversal and the separation between traversals is $2k+1$.
		
		\item In a type-$d$ discrete up-and-down path where $d$ is an even number in $[2k]$, the distance between consecutive stops along a traversal is $d$ and the separation between traversals is $2k+1-\frac{d}{2}$.
		
		\item We create type-$(2k+1)$ discrete up-and-down paths (zigzag paths) with a four step process, shown in Figure \ref{fig:DUAD}. As is done in \cite{carlsson2014continuous} to minimize stops, we first form a tessellation. Let $\mathbb{A}=\{(a,b)\in \mathbb{Z}^{2}|~\frac{ka+(k+1)b}{2k^{2}+2k+1}\in \mathbb{Z}\}$ and let $\mathbb{B}_{\mathbb{Z}}\big((a,b),k\big)$ be the set of integer points covered by $(a,b)$. Note that $\bigcup_{(a,b)\in \mathbb{A}}\mathbb{B}_{\mathbb{Z}}\big((a,b),k\big)=\mathbb{Z}^{2}$ and $\mathbb{B}_{\mathbb{Z}}\big((a_1,b_1),k\big)\cap \mathbb{B}_{\mathbb{Z}}\big((a_2,b_2),k\big)=\emptyset$ for any $(a_1,b_1)\neq (a_2,b_2)$ in $\mathbb{A}$. The stops in $\mathbb{A}$ along with their coverage define a tessellation of $\mathbb{Z}^{2}$ as shown in Figure \ref{fig:DUAD}(a). We choose all points $(a,b)\in \mathbb{A}$ such that $\mathbb{B}_{\mathbb{Z}}\big((a,b),k\big)\cap D\neq \emptyset$. Let $\mathbb{A}_{1}$ be the set of chosen stops. Next we determine the visit order for the stops in $\mathbb{A}_{1}$ using rotated traversals (Figure \ref{fig:DUAD}(b)). Traversal $i$ is the line segment connecting all stops $(a,b)$ in $\mathbb{A}_{1}$ satisfying $\frac{ka+(k+1)b}{2k^{2}+2k+1}=i$. Similar to Definition \ref{def:UAD}, we first connect these traversals in an up-and-down fashion which determines the order of stop connection. The tessellation may contain stops outside the grid; therefore, in the third step, we project those stops back to the grid. For each stop in $\mathbb{A}_{1}$ lying outside the rectangle $D$, we replace it with its projection onto the rectangle (Figure \ref{fig:DUAD}(c), black dots after projection are new stop locations). In the final step, the stops are connected using the grid lines according to the order defined in step 2 (Figure \ref{fig:DUAD}(d)).
	\end{itemize}
\end{definition}

We show in Corollaries \ref{mintour} and \ref{minstop} that the type-1 discrete up-and-down path and type-$(2k+1)$ discrete up-and-down path minimize path length $L$ and stop count $T$, respectively.

Propositions \ref{prop:feas-DUAD} and \ref{prop:cost-DUAD} establish the feasibility of discrete up-and-down path and compute its costs, respectively.
\begin{proposition}[Feasibility of discrete up-and-down path\label{prop:feas-DUAD}]
	For $d=1,2,4,\cdots,2k,2k+1$, a type-$d$ discrete up-and-down path covers all integer points in the rectangle $D$.
\end{proposition}

\begin{proposition}[Cost of discrete up-and-down path\label{prop:cost-DUAD}]
	The cost of discrete up-and-down path is summarized in Table \ref{table:2}.
\end{proposition}

\begin{table}[h]
	\begin{center}
		\newcommand{\tabincell}[2]{\begin{tabular}{@{}#1@{}}#2\end{tabular}}
		\begin{tabular}{ |c|c|c|c| }
			\hline
			Cost/Path type & type-$1$  & type-$2t$, $t=1,2,\cdots,k$ & type-$(2k+1)$ \\ \hline
			\tabincell{l}{~\\$L$\\~} & $\frac{N}{2k+1}+O(m)$  & $\frac{N}{2k+1-t}+O(m)$ & $\frac{(2k+1)N}{2k^{2}+2k+1}+O(km)$ \\ \hline
			\tabincell{l}{~\\$T$\\~} & $\frac{N}{2k+1}+O(m)$ & $\frac{N}{2t(2k+1-t)}+O(m)$ & $\frac{N}{2k^{2}+2k+1}+O(m)$ \\ \hline
		\end{tabular}
		\caption {Cost of discrete up-and-down path}\label{table:2}
	\end{center}
\end{table}
Similar to Proposition \ref{prop:cost-UAD}, the costs of type-1 and type-$2t$ discrete up-and-down paths are calculated based on the distance between consecutive stops and separation between traversals. For type-$(2k+1)$ discrete up-and-down path, the distance between (almost all) consecutive stops is $2k+1$ and almost every stop covers a set of distinct $2k^{2}+2k+1$ integer points. Thus, the stop count and path length are $\frac{N}{2k^{2}+2k+1}+O(m)$ and $\frac{(2k+1)N}{2k^{2}+2k+1}+O(km)$, respectively.
\begin{theorem}[\textbf{Tightness of trade-off inequality in the D-CPPG}\label{tightness-D}] Define the following piecewise-linear function $f_{UB-D}(\cdot)$ with turning points defined as
	\begin{equation}
		f_{UB-D}(d)=
		\begin{cases}
			2k+1& \text{if}~~d=1;\\
			d(2k+1-\frac{d}{2})& \text{if}~~d=2,4,\cdots,2k;\\
			2k^{2}+2k+1& \text{if}~~d\geq 2k+1.
		\end{cases}
	\end{equation}
	For any point $(L,T)$ satisfying
	\begin{equation}\label{eqn:30}
		N=(T-1)f_{UB-D}\Big(\frac{L}{T-1}\Big),
	\end{equation}
	there exists a feasible pair $(L^{'},T^{'})$ such that $L-L^{'}=O(km)$ and $T-T^{'}=O(km)$.
	
	Moreover, \eqref{eqn:30} is a polyline connecting
	\begin{equation}\label{eqn:polyline-UBD}
		\begin{split} &\Big(\frac{N}{f_{UB-D}(1)},\frac{N}{f_{UB-D}(1)}\Big)\rightarrow\Big(\frac{2N}{f_{UB-D}(2)},\frac{N}{f_{UB-D}(2)}\Big)\rightarrow\cdots\rightarrow\Big(\frac{2iN}{f_{UB-D}(2i)},\frac{N}{f_{UB-D}(2i)}\Big)\rightarrow\cdots\\
			&\rightarrow\Big(\frac{2kN}{f_{UB-D}(2k)},\frac{N}{f_{UB-D}(2k)}\Big)\rightarrow\Big(\frac{(2k+1)N}{f_{UB-D}(2k+1)},\frac{N}{f_{UB-D}(2k+1)}\Big)\rightarrow\Big(\infty,\frac{N}{f_{UB-D}(2k+1)}\Big).
		\end{split}
	\end{equation}
\end{theorem}

For each point $(L,T)$ on \eqref{eqn:polyline-UBD}, we construct a feasible mixed discrete up-and-down path that combines two different types of discrete up-and-down paths to cover $D_{int}$. We show that the costs of the mixed discrete up-and-down path are close to that of $(L,T)$. Therefore, \eqref{eqn:polyline-UBD} serves as an approximate upper bound of the feasible region in the D-CPPG.

\subsection{Optimization problem for the D-CPPG}
The following lemma shows the relative gap between polylines \eqref{eqn:polyline-D} and \eqref{eqn:polyline-UBD}, which are the lower bound and approximate upper bound of the feasible region of the D-CPPG, respectively.
\begin{lemma}\label{lem:gap-disc}
	In the D-CPPG, for any point $(L_1,T_1)$ on the lower bound polyline \eqref{eqn:polyline-D}, there exists $(L_2,T_2)$ on the approximate upper bound polyline \eqref{eqn:polyline-UBD} such that $\frac{L_2}{L_1}\leq \frac{11}{10},\frac{T_2}{T_1}\leq \frac{11}{10}$.
\end{lemma}

The following theorem follows from Lemma \ref{lem:gap-disc}, providing approximation ratio for optimization problem with linear objective.
\begin{theorem}[\textbf{Approximation ratio for the optimization problem}\label{approx-D}]
	If $C(L,T)=\alpha L+\beta T$ is a linear function and $m\geq n\geq O(\frac{k}{\varepsilon})$ where $\varepsilon\in (0,1)$, the mixed discrete up-and-down path provides a $\big(\frac{11}{10}+O(\varepsilon)\big)$-approximation solution.
\end{theorem}

\begin{corollary}[\textbf{Minimum path length in the D-CPPG}\label{mintour}]
	If $C(L,T)=C(L)$, the optimal solution is achieved with $L^{*}=\frac{N}{2k+1}+O(m)$.
\end{corollary}

\emph{Proof.}
From \eqref{eqn:polyline-D}, $L\geq \frac{N-(2k^{2}+2k+1)}{f_{LB-D}(1)}=\frac{N-(2k^{2}+2k+1)}{2k+1}$ for any covering path. A type-$1$ discrete up-and-down path achieves this bound.
\qed

In comparison with Corollary \ref{coro:1}, the $(2k+1)$ in the denominator is the maximum number of integer points covered by $F_{i+1}$ but not by $F_{i}$ given that the two stops are at distance 1. Locating a stop at each integer point of a traversal maximizes the separation between traversals, thus minimizing path length. Unlike in the C-CPPG, with a separation of $2k+1$ (rather than $2k$) we can derive a feasible path for $d=1$.

\begin{corollary}[\textbf{Minimum stop count in the D-CPPG}\label{minstop}]
	If $C(L,T)=C(T)$, the optimal solution is achieved with $T^{*}=\frac{N}{2k^{2}+2k+1}+O(m)$.
\end{corollary}
In Corollary \ref{coro:2}, the optimal solution is $T^{*}=\frac{N}{2k^{2}}+O(m)$. Here, we replace the denominator $2k^{2}$ with $(2k^{2}+2k+1)$, which is the number of integer points covered by a single stop in the D-CPPG.

\emph{Proof.}
From \eqref{eqn:to-disc} and \eqref{eqn:trade-off-disc}, $T-1\geq \frac{N-(2k^{2}+2k+1)}{2k^{2}+2k+1}$ for any feasible $T$, that is, $T\geq \frac{N}{2k^{2}+2k+1}$. A type-$(2k+1)$ discrete up-and-down path achieves this bound.
\qed

\subsection{Summary  of results}
In this subsection, we summarize the approximation results. Recall that $N = mn$ is the grid size and $k$ is the coverage radius. CPPG is trivially solvable when $k<1$ because a stop must be set at every grid vertex. Thus, we can assume that $k\geq 1$. As detailed in Figure \ref{fig:solapp}, CPPG can be solved by either solving C-CPPG or D-CPPG based on the value of $k$. For each of these settings we provide a feasible path where the parameters of the path are determined in polynomial time by solving a convex relaxation of the original optimization problem. Table \ref{tab:approx} summarizes our results for the three CPPG variants under different objective functions.

\begin{table}[h]
	\begin{center}
		\newcommand{\tabincell}[2]{\begin{tabular}{@{}#1@{}}#2\end{tabular}}
		\begin{tabular}{ |c|c|c|c| }
			\hline
			Objective & RC-CPPG & C-CPPG & D-CPPG \\ \hline
			$L$ & \tabincell{l}{\\~~~~~~~~$\frac{N}{2k}+O(m)$ \\ \\ $(1+\varepsilon)$-approximation\\~}  & \tabincell{l}{\\~~~~~~~~$\frac{N}{2k-\frac{1}{2}}+O(m)\sim \frac{N}{2k-1}+O(m)$~~~~~~~~~ \\ \\ ~~$(1+\frac{1}{4k-2})$-approximation ($\frac{11}{10}$ for $k\geq 3$)\\~} & \tabincell{l}{\\~~~~~~$\frac{N}{2k+1}+O(m)$ \\ \\ $(1+\varepsilon)$-approximation\\~}  \\ \hline
			$T$ & \tabincell{l}{\\~~~~~~~~$\frac{N}{2k^{2}}+O(m)$ \\ \\ $(1+\varepsilon)$-approximation\\~} & \tabincell{l}{\\~~~~~~~~$\frac{N}{2k^{2}}+O(m)$ \\ \\ $(1+\varepsilon)$-approximation\\~} & \tabincell{l}{\\~~~$\frac{N}{2k^{2}+2k+1}+O(m)$ \\ \\ $(1+\varepsilon)$-approximation\\~} \\ \hline
			$\alpha L+\beta T$ & \tabincell{l}{~\\$(1+\varepsilon)$-approximation\\~} & $(\frac{9}{8}+\varepsilon)$-approximation for $k\geq 3$ & $(\frac{11}{10}+\varepsilon)$-approximation \\ \hline
		\end{tabular}
		\caption {Summary of approximation results}
		\label{tab:approx}
	\end{center}
\end{table}

Observe that the only setting for which we do not have a $1+\varepsilon$ approximation is when minimizing $L$ for C-CPPG. Our approximation is weakest ($\frac{3}{2}$) for $k=1$ but strengthens as $k$ increases. For $k=3$, we obtain an $\frac{9}{8}$-approximation for CPPG for the general objective function and the approximation gets even stronger for larger values of $k$. As a result, in the worst case, we have a $\frac{3}{2}$-approximation for CPPG (when $k=1$) but for larger values of $k$ we are guaranteed much better results.

\section{Conclusion \label{conclusion}}
The core sub-problem in school bus routing is to select bus stops and a bus route connecting the stops such that no student is too far from a stop and the total bus route duration, including travel and stopping time, is minimized.  Motivated by the grid road structure of many American cities, we model the problem as one of obtaining a minimum cost covering path where the underlying network is a grid and distances are measured with the $l_1$ metric. Although the problem is known to be NP-hard on general graphs, we exploit the underlying grid structure to obtain strong approximations in polynomial time. Our solution approach is likely to be particularly useful as part of a decision support system where decision makers interactively build school bus routes by changing various input parameters.

We also feel that our results on complete unit grid graphs can become important building blocks for solution procedures on general grids. As long as the general grid can be constructed as the union of a few rectangular grids, our constructive approach can be used to find a solution that is unlikely to be too far from optimal. Our approach can also be used in the multi-vehicle setting that accounts for bus capacity. As long as we can divide the overall region into a union of rectangular grids (that relate to bus capacity), our constructive approach can be used to obtain a high quality solution for the capacitated multi-vehicle problem. In ongoing work, we are exploring ways to generalize the insights and results of this paper to address the many additional complications of school bus routing. We continue to work with the school district to provide solutions that are robust and easy to implement, and embed our results in to larger decision making frameworks for broader questions of school assignment.

\textbf{Acknowledgment}
	This project is supported by the National Science Foundation (CMMI-1727744) and a seed grant from the Office of Neighborhood and Community Relations at Northwestern University. The authors thank Superintendent Paul Goren, Assistant Superintendent Mary Brown, Transportation Coordinator Walter Doughty, and the staff of Evanston / Skokie Public School District 65.

\newpage
\setcounter{page}{1}
\small

\bibliographystyle{abbrvnat}
\bibliography{Literature}

\newpage
\normalsize

\section{Proof of Analytical Results in Section \ref{continuous}: Continuous CPPG}
\subsection{Proof of Lemma \ref{lem:piecewise}}
Since $g\big(\frac{X}{Y}\big)$ is a linear function of $\frac{X}{Y}$ when $\frac{X}{Y}\in [a_{i},a_{i+1})$, $Yg\big(\frac{X}{Y}\big)$ is a linear function of $X$ and $Y$ under the same condition. Therefore, given $\frac{X}{Y}\in [a_{i},a_{i+1})$ and constant $C>0$, $Yg\big(\frac{X}{Y}\big)=C$ is a line segment connecting $\big(\frac{a_{i}C}{b_{i}},\frac{C}{b_{i}}\big)$ and $\big(\frac{a_{i+1}C}{b_{i+1}},\frac{C}{b_{i+1}}\big)$. Summing over all cases of $i\in [n]$ , $Yg\big(\frac{X}{Y}\big)=C$ is equivalent to polyline \eqref{eqn:polyline-C}.

Furthermore, when $g(\cdot)$ is a piecewise-linear concave function, $g(\cdot)$ can be reformulated as the minimum of finite linear functions; i.e., $g(x)=\textrm{min}_{i}\{a_{i}x+b_{i} \}$. Therefore, $Yg\big(\frac{X}{Y}\big)=C$ is equivalent to $\textrm{min}_{i}\{a_{i}X+b_{i}Y\}=C$, which must be a piecewise-linear convex function.
\subsection{Proof of Theorem \ref{to-cont-C}}
We follow the same steps as in the proof of Theorem \ref{to-cont-R} until inequality \eqref{eqn:second-to-last},
\[N\leq |S_1|+\sum_{i=1}^{T-1}|S_{i+1}-S_{i}| \leq 2k^2+\sum_{i=1}^{T-1}f(d_i).\]
Since stops are located at integer points in the C-CPPG, all $d_{i}$ must be integers. From \eqref{eqn:fLB-C}, $f_{LB-C}(d_i)=f(d_i)$. Together with the concavity of $f_{LB-C}(\cdot)$ we have
\[N-2k^2\leq\sum_{i=1}^{T-1}f(d_i)=\sum_{i=1}^{T-1}f_{LB-C}(d_i) \leq (T-1)f_{LB-C}\Big(\frac{\sum_{i=1}^{T-1}d_i}{T-1}\Big)=(T-1)f_{LB-C}\Big(\frac{L}{T-1}\Big).\]

For the boundary of \eqref{eqn:trade-off-cont}, note that $f_{LB-C}(\cdot)$ is a piecewise-linear concave function. From Lemma \ref{lem:piecewise} and \eqref{eqn:boudary-LBC}, $(T-1)f_{LB-C}\big(\frac{L}{T-1}\big)=N-2k^{2}$ is equivalent to polyline \eqref{eqn:polyline-LBC}.
\subsection{Proof of Proposition \ref{prop:cost-MUAD}}
The mixed up-and-down path can be divided into three parts: the type-$d$ up-and-down path, the type-$(d+2)$ up-and-down path and the segment connecting these two paths. We estimate the cost of the type-$d$ and the type-$(d+2)$ paths based on Proposition \ref{prop:cost-UAD}.

From Proposition \ref{prop:cost-UAD}, the path length and stop count of the type-$d$ up-and-down path covering a $\left\lceil \gamma n\right\rceil \times m$ rectangle is at most $\big(\frac{m\left\lceil \gamma n\right\rceil}{2k-d/2}+3m\big)$ and $\big(\frac{\left\lceil\gamma n\right\rceil}{2k-d/2}+2\big)\big(\frac{m}{d}+2\big)$, respectively; the path length and stop count of the the type-$(d+2)$ up-and-down path covering a $\left\lceil (1-\gamma) n\right\rceil \times m$ rectangle is at most $\Big(\frac{m\left\lceil (1-\gamma) n\right\rceil}{2k-(d+2)/2}+3m\Big)$ and $\Big(\frac{\left\lceil (1-\gamma) n\right\rceil}{2k-(d+2)/2}+2\Big)\Big(\frac{m}{d+2}+2\Big)$, respectively.

Note that the length of the segment connecting these two paths is at most $2m$, the total path length of a type-$(d,\gamma)$ up-and-down path is at most
\begin{equation}\label{eqn:total-path}
	\begin{split}
		&\Big(\frac{m\left\lceil \gamma n\right\rceil}{2k-d/2}+3m\Big)+\Big(\frac{m\left\lceil (1-\gamma) n\right\rceil}{2k-(d+2)/2}+3m\Big)+2m\\
		& \leq \frac{m(\gamma n+1)}{2k-d/2}+\frac{m((1-\gamma) n+1)}{2k-(d+2)/2}+8m\\
		& =\frac{\gamma mn}{2k-d/2}+\frac{(1-\gamma) mn}{2k-(d+2)/2}+\frac{m}{2k-d/2}+\frac{m}{2k-(d+2)/2}+8m\\
		& \leq \frac{\gamma mn}{2k-d/2}+\frac{(1-\gamma) mn}{2k-(d+2)/2}+10m.
	\end{split}
\end{equation}
And the total stop count is at most
\begin{equation}\label{total-stop}
	\begin{split}
		& \Big(\frac{\left\lceil\gamma n\right\rceil}{2k-d/2}+2\Big)\Big(\frac{m}{d}+2\Big)+\Big(\frac{\left\lceil(1-\gamma)n\right\rceil}{2k-(d+2)/2}+2\Big)\Big(\frac{m}{d+2}+2\Big)\\
		& \leq \Big(\frac{\gamma n+1}{2k-d/2}+2\Big)\Big(\frac{m}{d}+2\Big)+\Big(\frac{(1-\gamma)n+1}{2k-(d+2)/2}+2\Big)\Big(\frac{m}{d+2}+2\Big)\\
		& \leq \Big(\frac{\gamma n}{2k-d/2}+3\Big)\Big(\frac{m}{d}+2\Big)+\Big(\frac{(1-\gamma)n}{2k-(d+2)/2}+3\Big)\Big(\frac{m}{d+2}+2\Big)\\
		& \leq \frac{\gamma mn}{d(2k-d/2)}+\frac{(1-\gamma) mn}{(d+2)(2k-(d+2)/2)}+\Big(\frac{3m}{d}+\frac{2\gamma n}{2k-d/2}+\frac{3m}{d+2}+\frac{2(1-\gamma)n}{2k-(d+2)/2}+12\Big)\\
		& \leq \frac{\gamma mn}{d(2k-d/2)}+\frac{(1-\gamma) mn}{(d+2)(2k-(d+2)/2)}+(3m+2n+3m+2n+12)\\
		& \leq \frac{\gamma mn}{d(2k-d/2)}+\frac{(1-\gamma) mn}{(d+2)(2k-(d+2)/2)}+10m+12.
	\end{split}
\end{equation}
\subsection{Proof of Theorem \ref{tight-cont-C}}
$f_{UB-C}(\cdot)$ is a concave function because of the concavity of $f(\cdot)$. From Lemma \ref{lem:piecewise}, equation \eqref{eqn:trade-off-C-eq} is equivalent to polyline \eqref{eqn:poly-UBC}. For any point $(L,T)$ on \eqref{eqn:poly-UBC}, we discuss two cases based on whether the point lies on the last segment of \eqref{eqn:poly-UBC}; i.e., $\big(\frac{2kN^{*}}{f(2k)},\frac{N^{*}}{f(2k)}\big)\rightarrow\big(\infty,\frac{N^{*}}{f(\infty)}\big)$.

\emph{Case 1:} If $(L,T)$ is not on $\big(\frac{2kN^{*}}{f(2k)},\frac{N^{*}}{f(2k)}\big)\rightarrow\big(\infty,\frac{N^{*}}{f(\infty)}\big)$, there exists $\gamma\in [0,1)$ and $t\in [k-1]$ such that
\begin{equation}
	L=\gamma \frac{2tN^{*}}{f(2t)}+(1-\gamma) \frac{(2t+2)N^{*}}{f(2t+2)},
\end{equation}
and
\begin{equation}
	T=\gamma \frac{N^{*}}{f(2t)}+(1-\gamma) \frac{N^{*}}{f(2t+2)}+1.
\end{equation}
From $f(d)=d(2k-\frac{d}{2})$ for $d\leq 2k$, we have
\begin{equation}\label{eqn:21}
	L=\gamma \frac{2tN^{*}}{f(2t)}+(1-\gamma) \frac{(2t+2)N^{*}}{f(2t+2)}=\frac{\gamma N^{*}}{2k-t}+\frac{(1-\gamma) N^{*}}{2k-(t+1)},
\end{equation}
and
\begin{equation}\label{eqn:22}
	T=\gamma \frac{N^{*}}{f(2t)}+(1-\gamma) \frac{N^{*}}{f(2t+2)}+1=\frac{\gamma N^{*}}{2t(2k-t)}+\frac{(1-\gamma) N^{*}}{\big(2t+2\big)\big(2k-(t+1)\big)}.
\end{equation}
Compare \eqref{eqn:21} to \eqref{eqn:total-path}, \eqref{eqn:22} to \eqref{total-stop}, together with the fact that $N=N^{*}\big(1+\frac{2k^{2}}{N-2k^{2}}\big)$, the costs of a type-$(2t,\gamma)$ mixed up-and-down path satisfy \eqref{eqn:23} and \eqref{eqn:24}.

\emph{Case 2:} If $(L,T)$ lies on $\big(\frac{2kN^{*}}{f(2k)},\frac{N^{*}}{f(2k)}\big)\rightarrow\big(\infty,\frac{N^{*}}{f(\infty)}\big)$, from $f(\infty)=f(2k)$, we have $L\geq \frac{2kN^{*}}{f(2k)}$ and $T=\frac{N^{*}}{f(2k)}$. From the analysis of \emph{Case 1}, the costs of a type-$2k$ up-and-down path satisfy \eqref{eqn:23} and \eqref{eqn:24}. Since we only increase $L$ and $T$ in \emph{Case 2}, inequalities \eqref{eqn:23} and \eqref{eqn:24} still hold.
\subsection{Proof of Lemma \ref{approx-C}}
Note that the turning points of \eqref{eqn:polyline-LBC} contains that of \eqref{eqn:poly-UBC}, it suffices to consider the extreme case where $(L_1, T_1)=\big(\frac{(2t+1)N^{*}}{f(2t+1)},\frac{N^{*}}{f(2t+1)}\big)$, a turning point of \eqref{eqn:polyline-LBC} but not of \eqref{eqn:poly-UBC}, where $t$ is an integer such that $0\leq t\leq k-1$. We choose $(L_2, T_2)$ based on the value of $t$.

\emph{Case 1:} if $t=0$, we take $(L_2, T_2)=\big(\frac{2N^{*}}{f(2)},\frac{N^{*}}{f(2)}\big)$. Then $\frac{L_2}{L_1}=\frac{2f(1)}{f(2)}=\frac{4k-1}{4k-2}\leq \frac{11}{10}$ for $k\geq 3$ ($\frac{3}{2}$ for $k=1$ and $\frac{7}{6}$ for $k=2$), $\frac{T_2}{T_1}=\frac{f(1)}{f(2)}<1$.

\emph{Case 2:} if $t>0$, let $(L_2, T_2)=\frac{2k-t}{(2k-t)+(2k-t-1)}\big(\frac{2tN^{*}}{f(2t)},\frac{N^{*}}{f(2t)}\big)+\frac{2k-t-1}{(2k-t)+(2k-t-1)}\big(\frac{(2t+2)N^{*}}{f(2t+2)},\frac{N^{*}}{f(2t+2)}\big)$. Clearly $(L_2, T_2)$ lies on the line segment connecting $\big(\frac{(2t)N^{*}}{f(2t)},\frac{N^{*}}{f(2t)}\big)$ and $\big(\frac{(2t+2)N^{*}}{f(2t+2)},\frac{N^{*}}{f(2t+2)}\big)$, therefore also on polyline \eqref{eqn:poly-UBC}. From $f(d)=d(2k-\frac{d}{2})$ for $d\leq 2k$, we have
\begin{equation}
	\begin{split}
		L_2 & =\frac{2k-t}{(2k-t)+(2k-t-1)}\frac{2tN^{*}}{f(2t)}+\frac{2k-t-1}{(2k-t)+(2k-t-1)}\frac{(2t+2)N^{*}}{f(2t+2)}\\
		& =\frac{2N^{*}}{(2k-t)+(2k-t-1)}\\
		& =\frac{(2t+1)N^{*}}{f(2t+1)}=L_1,
	\end{split}
\end{equation}
and
\begin{equation}
	\begin{split}
		T_2 & =\frac{2k-t}{(2k-t)+(2k-t-1)}\frac{N^{*}}{f(2t)}+\frac{2k-t-1}{(2k-t)+(2k-t-1)}\frac{N^{*}}{f(2t+2)}\\
		& =\frac{1}{2(2k-(2t+1)/2)}\Big(\frac{N^{*}}{2t}+\frac{N^{*}}{2t+2}\Big)\\
		& =\frac{1}{2(2k-(2t+1)/2)}\cdot \frac{2(2t+1)N^{*}}{2t(2t+2)}\\
		& =\frac{N^{*}}{(2t+1)(2k-(2t+1)/2)}\cdot  \frac{(2t+1)^{2}}{2t(2t+2)}\\
		& =\frac{N^{*}}{f(2t+1)}\cdot  \frac{(2t+1)^{2}}{2t(2t+2)}\\
		& =\big(1+\frac{1}{2t(2t+2)}\big)T_{1}\\
		& \leq \frac{9}{8}T_1.
	\end{split}
\end{equation}
In summary, $\frac{L_2}{L_1}\leq \frac{11}{10}$ for $k\geq 3$ ($\frac{3}{2}$ for $k=1$ and $\frac{7}{6}$ for $k=2$) and $\frac{T_2}{T_1}\leq \frac{9}{8}$.
\subsection{Proof of Theorem \ref{opt-C}}
For each point $(L_1,T_1)$ on \eqref{eqn:polyline-LBC}, from Lemma \ref{approx-C}, there exists $(L_2,T_2)$ on \eqref{eqn:poly-UBC} such that $\frac{L_2}{L_1}\leq \frac{11}{10}$ for $k\geq 3$ ($\frac{3}{2}$ for $k=1$ and $\frac{7}{6}$ for $k=2$), $\frac{T_2}{T_1}\leq \frac{9}{8}$. From Theorem \ref{tight-cont-C}, there exists a feasible pair $(L_3,T_3)$ such that $L_{3}-L_{2}\leq \frac{2k^{2}}{N-2k^{2}}L_{2}+10m$ and $T_{3}-T_{2}\leq \frac{2k^{2}}{N-2k^{2}}T_{2}+10m+12$.

When $m\geq n\geq \frac{100k}{\varepsilon}$, $\frac{L_{3}}{L_{2}}\leq 1+\frac{2k^{2}}{N-2k^{2}}+\frac{10m}{L_{2}}$. Since $(L_2,T_2)$ lies on \eqref{eqn:poly-UBC} and $L_2\geq \frac{2N^{*}}{f(2)}=\frac{N^{*}}{2k-1}$,
\begin{equation}
	\begin{split}
		& \frac{L_{3}}{L_{2}}\leq 1+\frac{2k^{2}}{N-2k^{2}}+\frac{10m}{L_{2}}\\
		& \leq 1+\frac{2k^{2}}{N-2k^{2}}+\frac{10(2k-1)m}{N^{*}}\\
		& \leq 1+\frac{2k^{2}}{10000k^{2}/\varepsilon^{2}}-2k^{2}+\frac{20km}{100km/\varepsilon-2k^{2}}\\
		& \leq 1+\frac{\varepsilon}{4}+\frac{\varepsilon}{4}\\
		& \leq 1+\frac{\varepsilon}{2}.
	\end{split}
\end{equation}
From \eqref{eqn:poly-UBC}, $T_2\geq \frac{N^{*}}{f(2k)}=\frac{N^{*}}{2k^{2}}$, we have,
\begin{equation}
	\begin{split}
		& \frac{T_{3}}{T_{2}}\leq 1+\frac{2k^{2}}{N-2k^{2}}+\frac{10m}{T_{2}}\\
		& \leq 1+\frac{2k^{2}}{10000k^{2}/\varepsilon^{2}-2k^{2}}+\frac{20k^{2}m}{N^{*}}\\
		& \leq 1+\frac{\varepsilon}{4}+\frac{20k^{2}m}{100km/\varepsilon-2k^{2}}\\
		& \leq 1+\frac{\varepsilon}{4}+\frac{\varepsilon}{4}\\
		& \leq 1+\frac{\varepsilon}{2}.
	\end{split}
\end{equation}
In summary, $\frac{L_{3}}{L_{1}}=\frac{L_{3}}{L_{2}}\cdot\frac{L_{2}}{L_{1}}\leq \frac{11}{10}(1+\frac{\varepsilon}{2})\leq \frac{11}{10}+\varepsilon$ for $k\geq 3$ ($\frac{3}{2}$ for $k=1$ and $\frac{7}{6}$ for $k=2$) and $\frac{T_{3}}{T_{1}}=\frac{T_{3}}{T_{2}}\cdot\frac{T_{2}}{T_{1}}\leq \frac{9}{8}(1+\frac{\varepsilon}{2})\leq \frac{9}{8}+\varepsilon$.
\section{Proof of Analytical Results in Section \ref{discrete}: Discrete CPPG}
\subsection{Proof of Theorem \ref{tradeoff-D}}\label{sec:B1}
Let $F_1-F_2-\cdots-F_T$ be a covering path where $\{F_i\}_{i=1}^{T}$ are the stops. Denote $d_i$ the distance between $F_i$ and $F_{i+1}$, and $L=\sum_{i=1}^{T-1}d_i$ the path length. Let $S_i$ be the set of integer points covered by $F_i$ and $f_{LB-D}(\cdot)$ be the lower bound function defined in \eqref{eqn:trade-off-disc}. Similar to \eqref{eqn:second-to-last} we have
\[N\leq |S_1|+\sum_{i=1}^{T-1}|S_{i+1}-S_{i}|.\]
It suffices to show $|S_{i+1}-S_i|\leq f_{LB-D}(d_i)$ and $f_{LB-D}(\cdot)$ is concave.

For the first claim, recall the definition in Lemma \ref{coverlemma} that $\mathbb{B}\big((a,b),k\big)=\{(x,y)~|~|x-a|+|y-b|\leq k\}$. Denote $\mathbb{B}_{\mathbb{Z}}\big((a,b),k\big)=\mathbb{B}\big((a,b),k\big)\cap\mathbb{Z}^{2}$ the set of integer points in $\mathbb{B}\big((a,b),k\big)$, similar to Lemma \ref{coverlemma}, we prove the following lemma.
\begin{lemma}\label{coverlemma-D}
	For any $k>0$ and $(p, q)\in \mathbb{Z}^{2}$,
	\[|\mathbb{B}_{\mathbb{Z}}\big((0,0),k\big)\cap \mathbb{B}_{\mathbb{Z}}\big((|p|+|q|,0),k\big)|\leq |\mathbb{B}_{\mathbb{Z}}\big((0,0),k\big)\cap \mathbb{B}_{\mathbb{Z}}\big((p,q),k\big)|. \]
\end{lemma}
Assume WLOG that $p\geq q\geq 0$. From inequality \eqref{eqn:subset} we know the following subset property of intersections
\begin{equation}
	\mathbb{B}_{\mathbb{Z}}\big((0,0),k\big)\cap \mathbb{B}_{\mathbb{Z}}\big((|p|+|q|,0),k\big)\subseteq \mathbb{B}_{\mathbb{Z}}\big((0,0),k\big)\cap \mathbb{B}_{\mathbb{Z}}\big((p,q),k\big).
\end{equation}
Therefore, Lemma \ref{coverlemma-D} is correct.

Lemma \ref{coverlemma-D} implies that $|S_{i+1}\cap S_{i}|\geq |\mathbb{B}_{\mathbb{Z}}\big((0,0),k\big)\cap \mathbb{B}_{\mathbb{Z}}\big((d_i,0),k\big)|$, hence,
\begin{equation}
	|S_{i+1}-S_i|=|S_{i+1}|-|S_{i+1}\cap S_{i}|\leq (2k^{2}+2k+1)-|\mathbb{B}_{\mathbb{Z}}\big((0,0),k\big)\cap \mathbb{B}_{\mathbb{Z}}\big((d_i,0),k\big)|.
\end{equation}
Next we show that
\begin{equation}\label{eqn:44}
	|\mathbb{B}_{\mathbb{Z}}\big((0,0),k\big)\cap \mathbb{B}_{\mathbb{Z}}\big((d,0),k\big)|=(2k^{2}+2k+1)-f_{LB-D}(d)
\end{equation}
for any integer $d$ and that $f_{LB-D}(\cdot)$ is concave.

Note that
\begin{equation}\label{eqn:39}
	|\mathbb{B}_{\mathbb{Z}}\big((0,0),k\big)\cap \mathbb{B}_{\mathbb{Z}}\big((d,0),k\big)|=|\{(x,y)|~|x|+|y|\leq k, |x-d|+|y|\leq k, x, y\in \mathbb{Z} \}|,
\end{equation}
we discuss two cases based on the value of $d$ to prove \eqref{eqn:44}.

\emph{Case 1:}
If $d\geq 2k+1$, the coverage area do not overlap; the right-hand side of \eqref{eqn:39} is 0.

\emph{Case 2:}
If $d\leq 2k$, for given $x$, the number of possible $y$ satisfying \eqref{eqn:39} is $2\cdot\text{min}\{k-|x|, k-|x-d|\}+1$. Note that both $|x|$ and $|x-d|$ are no more than $k$, the range of $x$ is $[d-k, k]$ and the total number of solutions to \eqref{eqn:39} is
\begin{equation}
	|\{(x,y)|~|x|+|y|\leq k, |x-d|+|y|\leq k, x, y\in \mathbb{Z} \}|=\sum_{x=d-k}^{k}(2\text{min}\{k-|x|, k-|x-d|\}+1).
\end{equation}
Now we prove the size of intersection satisfies
\begin{equation}\label{eqn:47}
	A=\sum_{x=d-k}^{k}(2\text{min}\{k-|x|, k-|x-d|\}+1)=(2k^{2}+2k+1)-f_{LB-D}(d).
\end{equation}
We define a threshold value that determines the minimum function in \eqref{eqn:47}. Note that $\text{min}\{k-|x|, k-|x-d|\}=k-|x|$ if and only if $|x|\geq |x-d|$; i.e., $x\geq \frac{d}{2}$; otherwise, if $x< \frac{d}{2}$, the minimum takes $k-|x-d|$. This threshold value $\frac{d}{2}$ is in $[d-k,k]$ since $d\leq 2k$. Since the threshold value $\frac{d}{2}$ is a integer or half-integer, we further separate \emph{Case 2} into two subcases based on the integrality of $\frac{d}{2}$.

\emph{Case 2.1:} $d$ is an odd number in $[2k]$. Because $\frac{d}{2}$ is a half-integer, we separate $A$ into two parts: $x$ in $[d-k,\frac{d-1}{2}]$ and $x$ in $[\frac{d+1}{2},k]$.
\begin{equation}\label{eqn:odd-cover}
	\begin{split}
		A&=\sum_{x=d-k}^{\frac{d-1}{2}}(2\text{min}\{k-|x|, k-|x-d|\}+1)+\sum_{x=\frac{d+1}{2}}^{k}(2\text{min}\{k-|x|, k-|x-d|\}+1)\\
		&=\sum_{x=d-k}^{\frac{d-1}{2}}\big(2(k-|x-d|)+1\big)+\sum_{x=\frac{d+1}{2}}^{k}\big(2(k-|x|)+1\big)\\
		&=\sum_{x=d-k}^{\frac{d-1}{2}}\big(2(k+x-d)+1\big)+\sum_{x=\frac{d+1}{2}}^{k}\big(2(k-x)+1\big)~~~(\textrm{break down absolute value})\\
		&=A_1+A_2.
	\end{split}
\end{equation}
$A_1$ counts for the summatation over $x$ in $[d-k,\frac{d-1}{2}]$ and $A_2$ over $x$ in $[\frac{d+1}{2},k]$.

Note that $A_1$ is the sum of an arithmetic sequence with first term is $2(k+(d-k)-d)+1=1$, last term is $2(k+\frac{d-1}{2}-d)+1=2k-d$ and common difference 2, we have
\begin{equation}\label{eqn:A1}
	\begin{split}
		A_1&=\frac{1}{2}\Big(1+(2k-d)\Big)\Big(\frac{(2k-d+1)-1}{2}+1\Big)\\
		&=\frac{(2k-d+1)^{2}}{4}.
	\end{split}
\end{equation}

Similarly, $A_2$ is the sum of an arithmetic sequence with first term is $2(k-\frac{d+1}{2})+1=2k-d$, last term is $2(k-k)+1=1$ and common difference -2, we have
\begin{equation}\label{eqn:A2}
	\begin{split}
		A_2&=\frac{1}{2}\Big((2k-d)+1\Big)\Big(\frac{(1-(2k-d)}{-2}+1\Big)\\
		&=\frac{(2k-d+1)^{2}}{4}.
	\end{split}
\end{equation}
Combining \eqref{eqn:A1} and \eqref{eqn:A2} we have
\begin{equation}
	A=A_1+A_2=\frac{(2k-d+1)^{2}}{2}.
\end{equation}
Note that $f_{LB-D}(d)=d(2k+1-\frac{d}{2})+\frac{1}{2}$ when $d$ is odd, we have
\[A+f_{LB-D}(d)=\frac{(2k-d+1)^{2}}{2}+d(2k+1-\frac{d}{2})+\frac{1}{2}=2k^{2}+2k+1.\]
Therefore, the size of intersection $A=(2k^{2}+2k+1)-f_{LB-D}(d)$ for odd $d$.

\emph{Case 2.2:} $d$ is an even number in $[2k]$. Because $\frac{d}{2}$ is an integer, we separate $A$ into two parts: $x$ in $[d-k,\frac{d}{2}]$ and $x$ in $[\frac{d+2}{2},k]$.
\begin{equation}\label{eqn:even-cover}
	\begin{split}
		A&=\sum_{x=d-k}^{\frac{d}{2}}(2\text{min}\{k-|x|, k-|x-d|\}+1)+\sum_{x=\frac{d+2}{2}}^{k}(2\text{min}\{k-|x|, k-|x-d|\}+1)\\
		&=\sum_{x=d-k}^{\frac{d}{2}}\big(2(k-|x-d|)+1\big)+\sum_{x=\frac{d+2}{2}}^{k}\big(2(k-|x|)+1\big)\\
		&=\sum_{x=d-k}^{\frac{d}{2}}\big(2(k+x-d)+1\big)+\sum_{x=\frac{d+2}{2}}^{k}\big(2(k-x)+1\big)~~~(\textrm{break down absolute value})\\
		&=A_3+A_4.
	\end{split}
\end{equation}
$A_3$ counts the summatation over $x$ in $[d-k,\frac{d}{2}]$ and $A_4$ over $x$ in $[\frac{d+2}{2},k]$.

$A_3$ is the sum of an arithmetic sequence with first term is $2(k+(d-k)-d)+1=1$, last term is $2(k+\frac{d}{2}-d)+1=2k-d+1$ and common difference 2, we have
\begin{equation}\label{eqn:A3}
	\begin{split}
		A_1&=\frac{1}{2}\Big(1+(2k-d+1)\Big)\Big(\frac{(2k-d+1)-1}{2}+1\Big)\\
		&=\frac{(2k-d+2)^{2}}{4}.
	\end{split}
\end{equation}

Similarly, $A_4$ is the sum of an arithmetic sequence with first term is $2(k-\frac{d+2}{2})+1=2k-d-1$, last term is $2(k-k)+1=1$ and common difference -2, we have
\begin{equation}\label{eqn:A4}
	\begin{split}
		A_2&=\frac{1}{2}\big((2k-d-1)+1\big)\Big(\frac{(1-(2k-d-1)}{-2}+1\Big)\\
		&=\frac{(2k-d)^{2}}{4}.
	\end{split}
\end{equation}
Combining \eqref{eqn:A3} and \eqref{eqn:A4}, we have
\begin{equation}
	A=A_3+A_4=\frac{(2k-d+2)^{2}}{4}+\frac{(2k-d)^{2}}{4}.
\end{equation}
Note that $f_{LB-D}(d)=d(2k+1-\frac{d}{2})$ when $d$ is even,
\[A+f_{LB-D}(d)=\frac{(2k-d+2)^{2}}{4}+\frac{(2k-d)^{2}}{4}+d(2k+1-\frac{d}{2})=2k^{2}+2k+1.\]
Therefore, the size of intersection $A=(2k^{2}+2k+1)-f_{LB-D}(d)$ for even $d$.

Now we prove the concavity of $f_{LB-D}(\cdot)$. Define $\Delta f_{LB-D}(d)=f_{LB-D}(d+1)-f_{LB-D}(d)$, it suffices to show $\Delta f_{LB-D}(d)$ is non-increasing for integer $d$. When $d$ is an even number in $[2k]$, $\Delta f_{LB-D}(d)=2k+1-d$; when $d$ is an odd number in $[2k]$, $\Delta f_{LB-D}(d)=2k-d$. Combining the fact that $\Delta f_{LB-D}(d)=0$ for all $d>2k$, $\Delta f_{LB-D}(d)$ is non-increasing and $f_{LB-D}(\cdot)$ is a piecewise-linear concave function.
\subsection{Proof of Proposition \ref{prop:feas-DUAD}}
We prove the feasibility of a discrete up-and-down path based on its type.

\emph{Case 1:} For type-$1$ discrete up-and-down path, each stop along a traversal covers itself together with its left and right $k$ columns of the grid. Since the separation between traversals is $2k+1$, all integer points in the grid are covered.

\emph{Case 2:} For type-$2t$ discrete up-and-down path, where $t\in[k]$, we follow the proof of Proposition \ref{prop:feas-UAD}. For any integer point $(x,y)$, assume that $(x,y)$ lies between traversals $i$ and $i+1$. Let $(x_h,y_h)$ be the highest stop on these two traversals with $y_h\leq y$ and $(x_l,y_l)$ be the lowest stop on these two traversals with $y_l\geq y$ (see Figure \ref{fig:UAD-cont} for illustration). From the alternating pattern of stop location, we can always pick $(x_h,y_h)$ and $(x_l,y_l)$ such that they are on different traversals. Since the separation between traversals $i$ and $i+1$ is at most $2k+1-\frac{d}{2}$ (equal to $2k+1-\frac{d}{2}$ except for the rightmost one), we have $|x_h-x_l|\leq 2k+1-\frac{d}{2}$. Also, recall the alternating pattern of stop locations on traversals $i$ and $i+1$, we have $|y_h-y_l|\leq \frac{d}{2}$.

Since $y_h\leq y\leq y_l$ and $x$ is always between $x_h$ and $x_l$, we have
\[||(x,y)-(x_h,y_h)||_{1}+||(x,y)-(x_l,y_l)||_{1}=||(x_h,y_h)-(x_l,y_l)||_{1}\leq (2k+1-\frac{d}{2})+\frac{d}{2}=2k+1.\]
This implies at least one of $||(x,y)-(x_l,y_l)||_{1}$ and $||(x,y)-(x_h,y_h)||_{1}$ is at most $k$. Thus, $(x_l,y_l)$ or $(x_h,y_h)$ covers $(x,y)$.

\emph{Case 3:} For type-$(2k+1)$ discrete up-and-down path, the coverage constraint is satisfied due to the tessellation property of coverage regions in the zigzag path.
\subsection{Proof of Proposition \ref{prop:cost-DUAD}}
For type-1 discrete up-and-down path, the separation between consecutive traversals is $2k+1$ and each traversal has $m$ stops. Therefore, both stop count and path length are $\frac{N}{2k+1}+O(m)$.

For type-$2t$ discrete up-and-down path, $t\in [k]$, the distance between consecutive traversals, distance between consecutive stops on one traversal are $2k+1-t$ and $2t$, respectively. Therefore, the stop count is $\frac{N}{2t(2k+1-t)}+O(m)$, and path length is $\frac{N}{2k+1-t}+O(m)$.

For type-$(2k+1)$ discrete up-and-down path, since the coverage regions of stops do not overlap, the stop count is $\frac{N}{2k^{2}+2k+1}+O(m)$. Note that distance between connected stops is $2k+1$, the path length is $\frac{(2k+1)N}{2k^{2}+2k+1}+O(km)$.
\subsection{Proof of Theorem \ref{tightness-D}}
Note that $f_{UB-D}(\cdot)$ is a concave piecewise-linear function, from Lemma \ref{lem:piecewise} we know that equation \eqref{eqn:30} is equivalent to polyline \eqref{eqn:polyline-UBD}. To illustrate the tightness result, we construct a feasible covering path for each point on \eqref{eqn:polyline-UBD} such that the cost is close enough to this point.

The turning points of the boundary of \eqref{eqn:30}  are: $T(d)=\frac{N}{f_{UB-D}(d)},L(d)=\frac{Nd}{f_{UB-D}(d)}$, where $d=1,2,4,\cdots,2k-2,2k,2k+1$. Therefore, the boundary point $(L,T)$ must lie on a line segment connecting two consecutive turning points with the following form
\begin{equation}
	L=\gamma L(d_1)+(1-\gamma)L(d_2), T=\gamma T(d_1)+(1-\gamma)T(d_2),
\end{equation}
where $d_1$ and $d_2$ are consecutive terms in the set $\{1,2,4,\cdots,2k-2,2k,2k+1\}$.

Besides, the path length and stop count of a type-$d$ discrete up-and-down path are $L(d)+O(km)$ and $T(d)+O(km)$, respectively. Let $(L^{'},T^{'})$ be the costs of a mixture of type-$d_1$ and type-$d_2$ discrete up-and-down path where the type-$d_1$ path covers a $m$ by $\gamma n$ grid and the type-$d_2$ path covers the other $m$ by $(1-\gamma) n$ grid. An additional path cost of order $O(m)$ connects the two discrete up-and-down paths. The stop count of this path is
\begin{equation}
	L^{'}=\gamma \big(L(d_1)+O(km)\big)+\big(1-\gamma\big)\big(L(d_2)+O(km)\big)=L+O(km),
\end{equation}
and the path length is
\begin{equation}
	T^{'}=\gamma \big(T(d_1)+O(km)\big)+\big(1-\gamma\big)\big(T(d_2)+O(km)\big)=T+O(km).
\end{equation}
To sum up, the gap between $(L^{'},T^{'})$ and $(L,T)$ is of order $O(km)$.
\subsection{Proof of Lemma \ref{lem:gap-disc}}
Note that $f_{UB-D}(d)=f_{LB-D}(d)$ for $d=1,2,4,\cdots,2k-2,2k,2k+1$, the turning points of \eqref{eqn:polyline-UBD} is a subset of that of \eqref{eqn:polyline-D}. Therefore, it suffices to consider $(L_1,T_1)$ as a turning point of \eqref{eqn:polyline-D} but not \eqref{eqn:polyline-UBD}; i.e., $(L_1,T_1)=\big(\frac{(2t+1)N}{f_{LB-D}(2t+1)},\frac{N}{f_{LB-D}(2t+1)}\big)$, where $t\in[k-1]$.

Let
\[L_2=\gamma \frac{2tN}{f_{UB-D}(2t)}+(1-\gamma)\frac{(2t+2)N}{f_{UB-D}(2t+2)} ,\]
and
\[T_2=\gamma \frac{N}{f_{UB-D}(2t)}+(1-\gamma)\frac{N}{f_{UB-D}(2t+2)} ,\]
where $\gamma=\frac{f_{UB-D}(2t)}{f_{UB-D}(2t)+f_{UB-D}(2t+2)}\in (0,1)$. Then $(L_2,T_2)$ lies on polyline \eqref{eqn:polyline-UBD}.

Moreover,
\begin{equation}
	\begin{split}
		L_2&=\gamma \frac{2tN}{f_{UB-D}(2t)}+(1-\gamma)\frac{(2t+2)N}{f_{UB-D}(2t+2)}\\
		&=\frac{f_{UB-D}(2t)}{f_{UB-D}(2t)+f_{UB-D}(2t+2)}\frac{2tN}{f_{UB-D}(2t)}+\frac{f_{UB-D}(2t+2)}{f_{UB-D}(2t)+f_{UB-D}(2t+2)}\frac{(2t+2)N}{f_{UB-D}(2t+2)}\\
		&=\frac{2(2t+1)N}{f_{UB-D}(2t)+f_{UB-D}(2t+2)}\\
		&=\frac{2f_{LB-D}(2t+1)}{f_{UB-D}(2t)+f_{UB-D}(2t+2)}L_1,
	\end{split}
\end{equation}
and
\begin{equation}
	\begin{split}
		T_2&=\gamma \frac{N}{f_{UB-D}(2t)}+(1-\gamma)\frac{N}{f_{UB-D}(2t+2)}\\
		&=\frac{f_{UB-D}(2t)}{f_{UB-D}(2t)+f_{UB-D}(2t+2)}\frac{N}{f_{UB-D}(2t)}+\frac{f_{UB-D}(2t+2)}{f_{UB-D}(2t)+f_{UB-D}(2t+2)}\frac{N}{f_{UB-D}(2t+2)}\\
		&=\frac{2N}{f_{UB-D}(2t)+f_{UB-D}(2t+2)}\\
		&=\frac{2f_{LB-D}(2t+1)}{f_{UB-D}(2t)+f_{UB-D}(2t+2)}T_1.
	\end{split}
\end{equation}
Note that
\begin{equation}
	\begin{split}
		&~~~~~f_{UB-D}(2t)+f_{UB-D}(2t+2)-2f_{LB-D}(2t+1)\\
		&=f_{LB-D}(2t)+f_{LB-D}(2t+2)-2f_{LB-D}(2t+1)\\
		&=2t(2k+1-t)+(2t+2)(2k-t)-2(2t+1)(2k+1-\frac{2t+1}{2})-1\\
		&=-2,
	\end{split}
\end{equation}
we have,
\begin{equation}\label{eqn:62}
	\begin{split}
		\frac{L_2}{L_1}=\frac{T_2}{T_1}&=\frac{2f_{LB-D}(2t+1)}{f_{UB-D}(2t)+f_{UB-D}(2t+2)}\\
		&=\frac{2f_{LB-D}(2t+1)}{2f_{LB-D}(2t+1)-2}\\
		&=1+\frac{1}{f_{LB-D}(2t+1)-1}.
	\end{split}
\end{equation}
From \eqref{eqn:trade-off-disc}, $f_{LB-D}(\cdot)$ is an increasing function and
\[f_{LB-D}(2t+1)\geq f_{LB-D}(3)=3(2k+1-\frac{3}{2})+\frac{1}{2}=6k-1\geq 11.\]
This, together with \eqref{eqn:62}, proves $\frac{L_2}{L_1}\leq \frac{11}{10},\frac{T_2}{T_1}\leq \frac{11}{10}$.

\subsection{Proof of Theorem \ref{approx-D}}
Similar to the proof of Theorem \ref{opt-C}, we first solve the optimization problem of minimizing $C(L,T)$ subject to $(T-1)f_{LB-D}(\frac{L}{T-1})=N-(2k^{2}+2k+1)$; i.e., minimizing $C(L,T)$ over polyline \eqref{eqn:polyline-D}. Let $(L^{*},T^{*})$ be the optimal solution, $\alpha L^{*}+\beta T^{*}$ is a lower bound for the optimal function value. From Lemma \ref{lem:gap-disc}, there exists $(L_{1},T_{1})$ on polyline \eqref{eqn:polyline-UBD} such that $\frac{L_{1}}{L^{*}}\leq \frac{11}{10}$ and $\frac{T_{1}}{T^{*}}\leq \frac{11}{10}$. Applying Theorem \ref{tightness-D}, there exists a feasible pair $(L^{'},T^{'})$ such that $L^{'}-L_1\leq O(km), T^{'}-T_1\leq O(km)$. When $m,n\geq O\big(\frac{k}{\varepsilon}\big)$, both $L^{*}$ and $T^{*}$ are at least $O\big(\frac{km}{\varepsilon}\big)$. Therefore,
\[\frac{L^{'}}{L^{*}}=\frac{L^{'}}{L_{1}}\cdot\frac{L_{1}}{L^{*}}\leq \big(1+O(\varepsilon)\big)\cdot\frac{11}{10}=\frac{11}{10}+O(\varepsilon). \]
Similar results holds for $\frac{L^{'}}{L^{*}}$ and $(L^{'},T^{'})$ is a $\big(\frac{11}{10}+O(\varepsilon)\big)$-approximation solution.

\end{document}